\documentclass[11pt]{article}
\usepackage{amsmath,amssymb,bm}
\def\N{\mathbb{N}}
\def\Z{\mathbb{Z}}
\def\C{\mathbb{C}}
\def\R{\mathbb{R}}
\def\L{\mathcal{L}}
\def\D{\mathcal{D}}
\def\G{\mathcal{G}}
\def\CB{\mathcal{CB}}
\def\a{\mathbf{a}}
\def\b{\mathbf{b}}

\def\m{\mathbf{m}}

\newtheorem{theorem}{\hspace*{\parindent}Theorem}
\newtheorem{lemma}{\hspace*{\parindent}Lemma}

\title{Representations of hypergeometric functions for arbitrary parameter values and their use}
\author{D.B.\:Karp$^{\rm a}$\footnote{Corresponding author. E-mail: D.\:Karp -- \emph{dimkrp@gmail.com}, J.L.\:L\'{o}pez --  \emph{jl.lopez@unavarra.es}}~~and J.L.\:L\'{o}pez$^{\rm b}$
\\[10pt]\small{\textit{$\phantom{1}^a$Far Eastern Federal University, Vladivostok, Russia}}\\\small{\textit{$\phantom{1}^b$Dpto. de Ingenier\'{\i}a Matem\'{a}tica e Inform\'{a}tica, Universidad P\'{u}blica de Navarra and INAMAT, Navarra, Spain}}}
\date{}
\begin{document}
\maketitle

\begin{center}
\parbox{12cm}{
\small\textbf{Abstract.}  Integral representations of hypergeometric functions proved to be a very useful tool for studying their properties. The purpose of this paper is twofold.  First, we extend the known representations to arbitrary values of the parameters and show that the extended representations can be interpreted as examples of regularizations of integrals containing Meijer's $G$ function. Second, we give new applications of both, known and extended representations. These include: inverse factorial series expansion for the Gauss type function, new information about zeros of the Bessel and Kummer type functions, connection with radial positive definite functions and generalizations of Luke's inequalities for the Kummer and Gauss type functions.}
\end{center}

\bigskip

Keywords: \emph{generalized hypergeometric function, Meijer's $G$ function, integral representation, radial positive definite function, inverse factorial series, Hadamard finite part}

\bigskip

MSC2010: 33C20, 33C60, 33F05, 42A82, 65D20

\bigskip

\section{Introduction}

Throughout the paper we will use the standard definition of the generalized hypergeometric function ${_{p}F_q}$ as the sum of the series
\begin{equation}\label{eq:pFqdefined}
{_{p}F_q}\left(\left.\!\!\begin{array}{c}\a \\ \b\end{array}\right|z\!\right)={_{p}F_q}\left(\a;\b;z\right)
=\sum\limits_{n=0}^{\infty}\frac{(a_1)_n(a_2)_n\cdots(a_{p})_n}{(b_1)_n(b_2)_n\cdots(b_q)_nn!}z^n
\end{equation}
if $p\le{q}$, $z\in\C$. If $p=q+1$ the above series only converges in the open unit disk and ${_{p}F_q}(z)$ is defined as its analytic continuation for $z\in\C\!\setminus\![1,\infty)$.  Here $(a)_n=\Gamma(a+n)/\Gamma(a)$ denotes the rising factorial (or Pochhammer's symbol) and $\a=(a_1,\ldots,a_p)$, $\b=(b_1,\ldots,b_q)$
are (generally complex) parameter vectors, such that $-b_j\notin\N_0$, $j=1,\ldots,q$.  This last restriction can be easily removed by dividing both sides of (\ref{eq:pFqdefined}) by $\prod_{k=1}^{q}\Gamma(b_k)$.  The resulting function (known as the regularized generalized hypergeometric function) is entire in $\b$. One useful tool in the study hypergeometric functions is their integral representations.  Probably, the earliest such representation is given by Euler's integral
$$
{_{2}F_1}\left(\sigma,a;b;-z\right)=\frac{\Gamma(b)}{\Gamma(a)\Gamma(b-a)}\int_{0}^{1}\frac{t^{a-1}(1-t)^{b-a-1}}{(1+zt)^{\sigma}}dt,
$$
that is finite for $z\in\C\setminus(-\infty,-1]$, $\Re(b-a)>0$ and $\Re(a)>0$.  This formula can be interpreted as the generalized Stieltjes transform of the
beta density $t^{a-1}(1-t)^{b-a-1}$.  In a recent paper \cite{Koornwinder}, Koornwinder generalized this formula and related
it to fractional integration formulas and transmutation operators. See also references in \cite{Koornwinder} for the history of the subject. Similar formulas with the generalized Stieltjes transform replaced by the Laplace and cosine Fourier transforms are valid for ${_{1}F_1}$ and ${_{0}F_1}$, respectively.
It seems surprising that for $p>1$  the generalized Stieltjes transform representation of ${_{p+1}F_p}$ (as well as the Laplace and cosine
Fourier transform representations for ${_{p}F_p}$ and ${_{p-1}F_p}$) has been only derived in 1994 by Kiryakova
in her book \cite[Chapter~4]{KiryakovaBook} and the article \cite{Kiryakova97} by the same author.  Her method of proof involves consecutive fractional
integrations and requires the restrictions $b_j>a_j>0$ on parameters in (\ref{eq:Frepr}).
We rediscovered similar representation using a different method in \cite{KSJAT} and utilized it to derive various inequalities and monotonicity results
for ${_{p+1}F_p}$.  Next, we relaxed the restrictions $b_j>a_j>0$ by demonstrating in \cite[Theorem~2]{KPJMAA} that, for an arbitrary complex $\sigma$, the representation
\begin{equation}\label{eq:Frepr}
{_{p+1}F_p}\left(\left.\!\!\begin{array}{c}\sigma,\a\\
\b\end{array}\right|-z\!\right)=\!\!\int\limits_{0}^{1}\frac{\rho(s)ds}{(1+sz)^{\sigma}}
\end{equation}
holds with a summable function $\rho$ and $|\arg(1+z)|<\pi$ if and only
if $\Re{a_i}>0$ for $i=1,\ldots,p$ and $\Re{\psi_p}>0$, where $\psi_p:=\sum_{j=1}^{p}(b_j-a_j)$. In the affirmative case
\begin{equation}\label{eq:rho-G}
\rho(s)=\frac{\Gamma(\b)}{\Gamma(\a)}G^{p,0}_{p,p}\left(s\left|\begin{array}{l}\!\!\b-1\!\!\\\!\!\a-1\!\!\end{array}\right.\right),
\end{equation}
where $G^{p,0}_{p,p}$ is Meijer's $G$-function defined in (\ref{eq:G-defined}) below. Further details about this function can be found in \cite[Section~12.3]{BealsWong}, \cite[Section~5.3]{HTF1}, \cite[Chapter~1]{KilSaig}, \cite[Section~8.2]{PBM3} and \cite[Section~16.17]{NIST}.  In (\ref{eq:rho-G}) we have used the abbreviated notation $\Gamma(\a)$ to denote the product $\prod_{i=1}^{p}\Gamma(a_i)$. This convention will also be used in the sequel. The sums like $\b+\alpha$ for a scalar $\alpha$ and inequalities like $\a>0$    will always be understood element-wise, i.e. $\b+\alpha=(b_1+\alpha,\ldots,b_p+\alpha)$ and $\a>0$ means $a_i>0$ for all elements of $\a$. Using term-by-term integration and some properties of the $G$ function, it is also straightforward to derive the following formulas to be used below \cite[(11), (12)]{KarpJMS}:
\begin{equation}\label{eq:qFqLaplace}
{_{p}F_p}\!\left(\left.\!\!\begin{array}{c}\a\\
\b\end{array}\right|-z\!\right)
=\frac{\Gamma(\b)}{\Gamma(\a)}\int\limits_{0}^{1}e^{-zt}G^{p,0}_{p,p}\left(t\left|\begin{array}{l}\!\!\b\!\!\\\!\!\a\!\!\end{array}\right.\right)
\frac{dt}{t}
\end{equation}
for $\Re(\a)>0$ and $\Re(\psi_p)>0$; and
\begin{equation}\label{eq:q-1Fqcosine}
{_{p-1}F_p}\left(\left.\!\!\begin{array}{c}\a\\\b\end{array}\right|-z\!\right)
=\frac{\Gamma(\b)}{\sqrt{\pi}\Gamma(\a)}\int\limits_{0}^{1}\cos(2\sqrt{zt})G^{p,0}_{p,p}\left(t\left|\begin{array}{l}\!\!\b\!\!\\\!\!\a,1/2\!\!\end{array}\right.\right)
\frac{dt}{t}
\end{equation}
for $\Re(\a)>0$  and $\Re(\sum_{j=1}^{p}b_j-\sum_{j=1}^{p-1}a_j)>1/2$.  Representations  (\ref{eq:Frepr}), (\ref{eq:qFqLaplace}) and (\ref{eq:q-1Fqcosine}) have been unified and generalized in \cite[Theorem~1]{KarpJMS}.

Although formulas (\ref{eq:Frepr}), (\ref{eq:qFqLaplace}) and (\ref{eq:q-1Fqcosine}) may be useful under the above restrictions on parameters in some contexts, most of the interesting applications appear for the values of parameters that ensure nonnegativity of the weight $\rho(s)$. The weakest known condition sufficient for the function in (\ref{eq:rho-G}) to be nonnegative is given by
\begin{equation}\label{eq:vge0}
v_{\a,\b}(t)=\sum\limits_{j=1}^{p}(t^{a_j}-t^{b_j})\geq0 ~~\text{for all}~t\in[0,1],
\end{equation}
as explained in \cite[Theorem~2]{KarpJMS}.  Further details regarding inequality (\ref{eq:vge0}) are collected in Property~9 found in the Appendix of this paper.
It follows from \cite[Theorem~3]{KPJMAA} combined with \cite[Lemma~2.1]{GrinIsm} (see also \cite[Theorem~6.4]{Dufresne2010}) that conditions (\ref{eq:vge0})  and $\a>0$
are sufficient for the representation
\begin{equation}\label{eq:Frepr1}
{_{p+1}F_p}\left(\left.\!\!\begin{array}{c}\sigma,\a\\
\b\end{array}\right|-z\!\right)
=\!\!\int\limits_{0}^{1}\frac{\mu(ds)}{(1+sz)^{\sigma}}
\end{equation}
to hold with a nonnegative measure $\mu$ supported on $[0,1]$. Hence, in contrast to (\ref{eq:Frepr}), representation (\ref{eq:Frepr1}) may hold when the parametric excess $\psi_p=0$, since this condition is compatible  with (\ref{eq:vge0}) (in fact necessary for it).  This shows that the measure $\mu(ds)$ is not absolutely continuous with respect to the Lebesgue measure if $\psi_p=0$.  We have demonstrated in \cite[p.\:353]{KPJMAA} that for  $p=1,2$ this measure has an atom at $s=1$ and an absolutely continuous part (vanishing for $p=1$).  The same result was discovered two years earlier by Dufresne in \cite[Theorem~6.2]{Dufresne2010} in a probabilistic context.  The first aim of this paper is to generalize this result to arbitrary $p\ge1$ and supply explicit expressions for both the atom and the absolutely continuous part.  This is done in Section~2, which further studies the limit of the measure $\rho(s)ds$, with $\rho(s)$ given in (\ref{eq:rho-G}), as $\min(\a)\to0$, and representations of the Kummer and Gauss type functions for $\psi_p$ equal to a negative integer.  Such representations can be easily derived with the help of some subtle characteristics of $G$ function (due to N{\o}rlund) outlined in Property~7 in the Appendix.

Next, we present three new applications of the integral representations (\ref{eq:Frepr}), (\ref{eq:qFqLaplace}) and (\ref{eq:q-1Fqcosine}) and their limiting cases in Section~3.    They are: convergent inverse factorial series expansion for ${_{p+1}F_p}$ for general complex parameters (subsection~3.1), new information about zeros of the Bessel and Kummer type functions and inequalities for the former (subsection~3.2) and an investigation of hypergeometric functions as radial positive definite functions (subsection~3.3).

Section~4 is devoted to the case of unrestricted complex parameters.  Straightforward decomposition of the series (\ref{eq:pFqdefined}) combined with integral representations for the remainders leads to representations of generalized hypergeometric functions as sums of Taylor polynomials and generalized Stieltjes, Laplace or cosine Fourier transforms of a complex density, see Theorem~\ref{th:pFqdecompose}. An important feature of these representations is that for any real parameters $\a$ and $\b$ taking sufficiently large degree of the Taylor polynomial makes this density nonnegative.  We show further that such decomposed representations can be seen as examples of regularizations of the divergent integrals containing Meijer's $G$ function (\ref{eq:rho-G}). The regularization theory is  developed in subsection~4.2, where the relation to the Hadamard finite part integrals is also observed.   Finally, subsection~4.3 contains an application of the decomposition theorem to the extension of Luke's inequalities to arbitrary real parameters for the Kummer and Gauss type functions. The properties of Meijer's $G$ function used throughout the paper are collected in the Appendix.

\section{The parametric excess is a non-positive integer}

\subsection{Limits as $\psi_p\to0$ and $\min(\a)\to0$}

In this section we assume that the parameter vectors $\a$ and $\b$ are real.
Denote the unbounded closed set in $\R^{2p}$ defined by inequalities (\ref{eq:vge0}) by $\D$.  It follows from \cite[(19)]{KPCMFT} that the boundary of $\D$ contains points of the  hyperplane $\psi_p=\sum_{i=1}^{p}(b_i-a_i)=0$.  The expression $\psi_p\to0$ will mean that $\psi_p$ is approaching the points of $\partial{\D}$ belonging to the hyperplane $\psi_p=0$ along any curve lying entirely in $\D$. The next lemma is elementary and probably well-known. We found it easier, however, to give a proof than to locate one in the literature.

\begin{lemma}\label{lm:limit}
Suppose that $f(t)$, with $f(0)=0$, is continuous on $[0,1)$ and absolutely integrable on $(0,1)$. Then
$$
\lim\limits_{\beta\downarrow0}\frac{1}{\Gamma(\beta)}\int_0^1t^{\beta-1}f(t)dt=0.
$$
\end{lemma}
\textbf{Proof.} Take an arbitrary $\varepsilon>0$. Since $f(t)$ is continuous on $[0,1)$ and $f(0)=0$ there exists $\lambda\in(0,1)$ such that $|f(\eta)|<\varepsilon/4$ for all $\eta\in[0,\lambda]$.  Further, for this $\lambda$, and according to the mean value theorem we have:
$$
\left|\frac{1}{\Gamma(\beta)}\int_0^\lambda t^{\beta-1}f(t)dt\right|
=\left|\frac{f(\eta)}{\Gamma(\beta)}\int_0^\lambda t^{\beta-1}dt\right|
\leq\frac{|f(\eta)|}{\Gamma(\beta)}\int_0^1t^{\beta-1}dt
=\frac{|f(\eta)|}{\Gamma(\beta+1)}<2|f(\eta)|<\varepsilon/2,
$$
where we have used the fact that $1/2<\Gamma(\beta+1)$ for $\beta\in(0,1)$.  The above estimate is independent of $\beta$.  Since $[\Gamma(0)]^{-1}=0$ we can choose $\delta$ such that for all $0<\beta<\delta$:
$$
\left|\frac{1}{\Gamma(\beta)}\int_{\lambda}^1t^{\beta-1}f(t)dt\right|\leq\frac{1}{\Gamma(\beta)}
\int_{\lambda}^1|f(t)|\frac{dt}{t}<\frac{\varepsilon}{2}.
$$
Hence, for all $0<\beta<\delta$ we get $|[\Gamma(\beta)]^{-1}\int_0^1t^{\beta-1}f(t)dt|<\varepsilon$ which completes the proof.~~$\square$

The following theorem extends \cite[Theorem~6.2]{Dufresne2010} and \cite[Remark on p.354]{KPJMAA}.
\begin{theorem}\label{th:zerobalanced}
The family of the probability measures
$$
\rho(ds)=\frac{\Gamma(\b)}{\Gamma(\a)}G^{p,0}_{p,p}\!\left(\!s~\vline\begin{array}{l}\b\\\a\end{array}\!\!\right)\!\!\frac{ds}{s},
$$
supported on $[0,1]$, converges weakly to the measure
$$
\frac{\Gamma(\b^*)}{\Gamma(\a^*)}\left\{\delta_{1}+G^{p,0}_{p,p}\!\left(\!s~\vline\begin{array}{l}\b^*\\\a^*\end{array}\!\!\right)\!\!\frac{ds}{s}\right\}
~~\text{as}~ \psi_p\to0~\text{in}~\D,
$$
where $\delta_1$ denotes the unit mass concentrated at the point $s=1$, and $(\a^*,\b^*)$ is a point on the hyperplane $\sum_{i=1}^{p}(a_i^*-b_i^*)=0$ such that $\a^*=\lim_{\psi_p\downarrow0}\a$, $\b^*=\lim_{\psi_p\downarrow0}\b$\emph{;} the $G^{p,0}_{p,p}$ function in the last formula is given by the integral \emph{(\ref{eq:G-defined})} with $\L=\L_{-}$ and can be computed by expansion \emph{(\ref{eq:Norlund1})}.
\end{theorem}
\textbf{Proof.} According to the definition of weak convergence of measures \cite[Section~10.3]{Cohn}, we need to show that, for any continuous function $\varphi(s)$ on $[0,1]$,
$$
\lim\limits_{\psi_p\to0}\int\limits_{0}^{1}\varphi(s)\rho(ds)=
\frac{\Gamma(\b^*)}{\Gamma(\a^*)}\left\{\varphi(1)+\int_0^1\!\varphi(s)G^{p,0}_{p,p}\!\left(\!s~\vline\begin{array}{l}\b^*\\\a^*\end{array}\!\!\right)\!\frac{ds}{s}\right\}.
$$
Rewrite (\ref{eq:Norlund}) as
\begin{equation}\label{eq:Gp0ppform}
G^{p,0}_{p,p}\!\left(\!z~\vline\begin{array}{l}\b
\\\a\end{array}\!\!\right)=\frac{z^{a_k}(1-z)^{\psi_p-1}}{\Gamma(\psi_p)}\left[1+\Gamma(\psi_p)(1-z)G_k(\a,\b;z)\right],
\end{equation}
where $G_k(\a,\b;z)$ is regular around $z=1$ and around $\psi_p=0$. Set $\phi(t)=\varphi(1-t)$ and $\tilde{\phi}(t)=\phi(t)-\phi(0)$, so that $\tilde{\phi}(t)$ is continuous on $[0,1]$ and $\tilde{\phi}(0)=0$. Next, compute:
\begin{multline*}
\lim\limits_{\psi_p\to0}\int\limits_{0}^{1}\varphi(s)\rho(ds)=\frac{\Gamma(\b^*)}{\Gamma(\a^*)}\lim\limits_{\psi_p\to0}
\int\limits_{0}^{1}\frac{s^{a_k-1}(1-s)^{\psi_p-1}}{\Gamma(\psi_p)}\left[1+\Gamma(\psi_p)(1-s)G_k(\a,\b;s)\right]\varphi(s)ds
\\
=\frac{\Gamma(\b^*)}{\Gamma(\a^*)}\lim\limits_{\psi_p\to0}
\int\limits_{0}^{1}\frac{(1-t)^{a_k-1}t^{\psi_p-1}}{\Gamma(\psi_p)}\left[1+\Gamma(\psi_p)tG_k(\a,\b;1-t)\right]\phi(t)dt
\\
=\frac{\Gamma(\b^*)}{\Gamma(\a^*)}\biggl\{\lim\limits_{\psi_p\to0}
\int\limits_{0}^{1}\frac{(1-t)^{a_k-1}t^{\psi_p-1}}{\Gamma(\psi_p)}(\phi(0)+\tilde{\phi}(t))dt
+ \lim\limits_{\psi_p\to0}\int\limits_{0}^{1}(1-t)^{a_k-1}t^{\psi_p}G_k(\a,\b;1-t)\phi(t)dt\biggr\}
\\
=\frac{\Gamma(\b^*)}{\Gamma(\a^*)}\biggl\{\lim\limits_{\psi_p\to0}
\frac{\phi(0)}{\Gamma(\psi_p)}\int\limits_{0}^{1}(1-t)^{a_k-1}t^{\psi_p-1}dt+
\lim\limits_{\psi_p\to0}\frac{1}{\Gamma(\psi_p)}\int\limits_{0}^{1}(1-t)^{a_k-1}t^{\psi_p-1}\tilde{\phi}(t)dt
\\
+\lim\limits_{\psi_p\to0}\int\limits_{0}^{1}(1-t)^{a_k-1}t^{\psi_p}G_k(\a,\b;1-t)\phi(t)dt\biggr\}
\\
=\frac{\Gamma(\b^*)}{\Gamma(\a^*)}\biggl\{\lim\limits_{\psi_p\to0}
\frac{\phi(0)}{\Gamma(\psi_p)}\frac{\Gamma(a_k)\Gamma(\psi_p)}{\Gamma(a_k+\psi_p)}+
0+\int\limits_{0}^{1}\phi(t)\lim\limits_{\psi_p\to0}(1-t)^{a_k-1}G_k(\a,\b;1-t)dt\biggr\}
\\
=\frac{\Gamma(\b^*)}{\Gamma(\a^*)}\biggl\{
\phi(0)+\int\limits_{0}^{1}\phi(t)\lim\limits_{\psi_p\to0}(1-t)^{a_k-1}G_k(\a,\b;1-t)dt\biggr\}
\\
=\frac{\Gamma(\b^*)}{\Gamma(\a^*)}\biggl\{\varphi(1)+\int\limits_{0}^{1}\frac{\varphi(s)}{s}\lim\limits_{\psi_p\to0}s^{a_k}G_k(\a,\b;s)ds\biggr\}.
\end{multline*}
Further, from (\ref{eq:Gp0ppform}):
$$
\lim\limits_{\psi_p\to0}s^{a_k}G_k(\a,\b;s)=\lim\limits_{\psi_p\to0}\frac{1}{(1-s)^{\psi_p}}
\left[G^{p,0}_{p,p}\!\left(\!s~\vline\begin{array}{l}\b\\\a\end{array}\!\!\right)-
\frac{s^{a_k}(1-s)^{\psi_p-1}}{\Gamma(\psi_p)}\right]
=G^{p,0}_{p,p}\!\left(\!s~\vline\begin{array}{l}\b^*\\\a^*\end{array}\!\!\right).
$$
We applied Lemma~\ref{lm:limit} in the fifth equality of the above chain.$\hfill\square$

\textbf{Remark.}  The above theorem has been extended in our recent paper \cite[Theorem~3]{KPFOX} to Fox's $H$ function.
However, the above proof is different and more direct than the one given in \cite{KPFOX}.  In view of this fact and for the sake of completeness,
we decided to present a full proof here.

Setting $\varphi(s)=(1+zs)^{-\sigma}$ in the above theorem, we get the representation
$$
{_{p+1}F_p}\left(\left.\!\!\begin{array}{c}\sigma,\a\\
\b\end{array}\right|-z\!\right)=\frac{\Gamma(\b)}{\Gamma(\a)}\biggl\{\frac{1}{(1+z)^{\sigma}}
+\int\nolimits_{0}^{1}\frac{ds}{s(1+sz)^{\sigma}}G^{p,0}_{p,p}\!\left(\!s~\vline\begin{array}{l}\b\\\a\end{array}\!\!\right)\biggr\}
$$
valid for $\psi_p=\sum(b_i-a_i)=0$.  The restriction (\ref{eq:vge0}) is removed by analytic continuation.
All we need is the condition $\Re(\a)>0$ for the above integral to converge.  Similarly, setting $\varphi(s)=e^{-zs}$, we get (again for $\psi_p=0$):
$$
{_{p}F_p}\left(\left.\!\!\begin{array}{c}\a\\\b\end{array}\right|-z\!\right)
=\frac{\Gamma(\b)}{\Gamma(\a)}\biggl\{e^{-z}
+\int\nolimits_{0}^{1}e^{-zs}G^{p,0}_{p,p}\!\left(\!s~\vline\begin{array}{l}\b\\\a\end{array}\!\!\right)\frac{ds}{s}\biggr\}.
$$
Finally if $\varphi(s)=\cos(2\sqrt{zs})$ and $\psi_p=1/2$, then
\begin{equation*}
{_{p-1}F_p}\left(\left.\!\!\begin{array}{c}\a\\\b\end{array}\right|-z\!\right)
=\frac{\Gamma(\b)}{\sqrt{\pi}\Gamma(\a)}\biggl\{\cos(2\sqrt{z})
+\int\nolimits_{0}^{1}\cos(2\sqrt{zs})G^{p,0}_{p,p}\left(s\left|\begin{array}{l}\!\!\b\!\!\\\!\!\a,1/2\!\!\end{array}\right.\right)\frac{ds}{s}\biggr\}.
\end{equation*}
Note that the above representations are particular, $\psi_p=0$, cases of Theorem~\ref{th:repr-psiinteger} below.

\begin{theorem}\label{th:agoestozero}
Set $a=\min(a_1,a_2,\ldots,a_p)$.  The family of probability measures
$$
\rho(ds)=\frac{\Gamma(\b)}{\Gamma(\a)}G^{p,0}_{p,p}\!\left(\!s~\vline\begin{array}{l}\b\\\a\end{array}\!\!\right)\!\!\frac{ds}{s},
$$
supported on $[0,1]$, converges weakly to the Dirac measure $\delta_0$ \emph{(}unit mass at zero\emph{)} as $a\to0$ staying in $\D$.
\end{theorem}
\textbf{Proof.}  Indeed, (\ref{eq:GMellin}) shows that $\rho$ is indeed a probability measure for parameters in $\D$. According to the definition of weak convergence \cite[Section~10.3]{Cohn} we need to show that for any continuous function $\phi(s)$ on $[0,1]$
$$
\lim\limits_{a\to0}\int\limits_{0}^{1}\phi(s)\rho(ds)-\phi(0)=\lim\limits_{a\to0}\int\limits_{0}^{1}\tilde{\phi}(s)\rho(ds)=0,
$$
where $\tilde{\phi}(s)=\phi(s)-\phi(0)$.  Choose an arbitrary $\delta>0$. We will prove that there exists $\lambda>0$ such that for all $0<a<\lambda$,
\begin{equation}\label{eq:delta-est1}
\left|\int\limits_{0}^{1}\tilde{\phi}(s)\rho(ds)\right|<\delta.
\end{equation}
Since  $\tilde{\phi}(s)$ is continuous on $[0,1)$ with $\tilde{\phi}(0)=0$, there exits $\varepsilon>0$ such that $|\tilde{\phi}(\eta)|<\delta/2$ for all $\eta\in[0,\varepsilon]$.  Further, for this value of $\varepsilon$  the mean value theorem yields:
$$
\left|\int\limits_{0}^{\varepsilon}\rho(s)\tilde{\phi}(s)ds\right|=\left|\tilde{\phi}(\eta)\int\limits_{0}^{\varepsilon}\rho(s)ds\right|
\leq|\tilde{\phi}(\eta)|\int\limits_{0}^{1}\rho(s)ds<\delta/2.
$$
The above estimate is independent of $a$.  Now choose $\lambda$ such that for all $0<a<\lambda$:
$$
\left|\int\nolimits_{\varepsilon}^{1}\tilde{\phi}(s)\rho(ds)\right|
=\frac{\Gamma(\b)}{\Gamma(\a)}
\left|\int\nolimits_{\varepsilon}^{1}G^{p,0}_{p,p}\!\left(\!s~\vline\begin{array}{l}\b\\\a\end{array}\!\!\right)\tilde{\phi}(s)\frac{ds}{s}\right|<\frac{\delta}{2}.
$$
This is possible because $[\Gamma(\a)]^{-1}\to0$ as $a\to0$ and the integrand is bounded on $s\in[\varepsilon,1]$ uniformly in $a$.
Hence, for all $0<a<\lambda$ we get (\ref{eq:delta-est1}), which completes the proof.~~$\hfill\square$

\subsection{The parametric excess $\psi_p$ is a negative integer}

In this subsection we explore the consequences of N{\o}rlund's formula (\ref{eq:GMellinNorlund}) valid for
non-positive integer values of $\psi_p$. We will use the notation $\a_{[k]}=(a_1,\ldots,a_{k-1},a_{k+1},\ldots,a_p)$.
New representations derived from this formula are presented in the next theorem.

\begin{theorem}\label{th:repr-psiinteger}
Suppose $-\psi_p=m\in\N_0$ and $\a>0$. Then
\begin{multline}\label{eq:p+1Fp-negintpsi}
\frac{\Gamma(\a)}{\Gamma(\b)}{_{p+1}F_p}\left(\left.\!\!\begin{array}{c}\sigma,\a\\
\b\end{array}\right|-z\!\right)=\frac{\Gamma(a_k)}{(1+z)^{\sigma}}\sum\limits_{j=0}^{m}\frac{g_j(\a_{[k]};\b)}{\Gamma(a_k+j-m)}
{_{2}F_1}\left(\left.\!\!\begin{array}{c}\sigma,j-m\\
a_k+j-m\end{array}\right|\frac{z}{1+z}\!\right)
\\
+\int\limits_{0}^{1}(1+zs)^{-\sigma}G^{p,0}_{p,p}\!\left(\!s~\vline\begin{array}{l}\b-1\\\a-1\end{array}\!\!\right)ds
\end{multline}
for $z\in\C\!\setminus\!(-\infty,-1]$ and arbitrary complex $\sigma$; for all complex $z$,
\begin{multline}\label{eq:pFp-negintpsi}
\frac{\Gamma(\a)}{\Gamma(\b)}{_{p}F_p}\left(\left.\!\!\begin{array}{c}\a\\
\b\end{array}\right|-z\!\right)=e^{-z}\Gamma(a_k)\sum\limits_{j=0}^{m}\frac{g_j(\a_{[k]};\b)}{\Gamma(a_k+j-m)}
{_{1}F_1}\left(\left.\!\!\begin{array}{c}j-m\\
a_k+j-m\end{array}\right|z\!\right)
\\
+\int\limits_{0}^{1}e^{-zs}G^{p,0}_{p,p}\!\left(\!s~\vline\begin{array}{l}\b-1\\\a-1\end{array}\!\!\right)ds.
\end{multline}
If $\psi_p=\sum_{j=1}^{p}b_j-\sum_{j=1}^{p-1}a_j=-m+1/2$, $m\in\N_0$, then
\begin{multline}\label{eq:p-1Fp-negintpsi}
\frac{\sqrt{\pi}\Gamma(\a)}{\Gamma(\b)}{_{p-1}F_p}\left(\left.\!\!\begin{array}{c}\a\\
\b\end{array}\right|-z\!\right)=\sum\limits_{j=0}^{m}(-1)^j(1/2)_jg_{m-j}(\a;\b)
{_{0}F_1}\left(\left.\!\!\begin{array}{c}-\\
1/2-j\end{array}\right|-z\!\right)
\\
+\int\limits_{0}^{1}\cos\left(2\sqrt{zs}\right)G^{p,0}_{p,p}\!\left(\!s~\vline\begin{array}{l}\b-1\\\a-1,-1/2\end{array}\!\!\right)ds,
\end{multline}
for all complex $z$. Formulas \emph{(\ref{eq:p+1Fp-negintpsi})} and \emph{(\ref{eq:pFp-negintpsi})} are valid for each $k=1,\ldots,p$. The coefficients $g_j(\a_{[k]};\b)$ are defined by the recurrence \emph{(\ref{eq:Norlundcoeff})} and the connection formula \emph{(\ref{eq:Norlundconnection})}, or explicitly in \emph{(\ref{eq:Norlund-explicit})}.
\end{theorem}
\textbf{Proof.}  For the proof, substitute the power series expansions of $(1+zs)^{-\sigma}$, $e^{-zs}$ and $\cos\left(2\sqrt{zs}\right)$ into (\ref{eq:p+1Fp-negintpsi}), (\ref{eq:pFp-negintpsi}) and (\ref{eq:p-1Fp-negintpsi}), respectively and integrate term by term using (\ref{eq:GMellinNorlund}).  Then apply Pfaff's transformation \cite[formula~2.2.6]{AAR} to the the resulting $_2F_1$ in (\ref{eq:p+1Fp-negintpsi}) and Kummer's transformation \cite[formula~4.1.11]{AAR} to the resulting $_1F_1$ in (\ref{eq:pFp-negintpsi}).~~$\hfill\square$

\noindent\textbf{Remark.} Note that the functions $_2F_1$ and $_1F_1$ in (\ref{eq:p+1Fp-negintpsi}) and (\ref{eq:pFp-negintpsi}), respectively, are finite sums. Furthermore, the ${}_0F_1$ in (\ref{eq:p-1Fp-negintpsi}) can be expressed as $\cos(2\sqrt{z})$ times a combination of Lommel polynomials.

For $p=2$ and  $-\psi_p\notin\N_0$ we have \cite[8.4.49.22]{PBM3}:
\begin{equation}\label{eq:Meijer2reduced}
G^{2,0}_{2,2}\!\left(\!t~\vline\begin{array}{l}b_1,b_2
\\a_1,a_2\end{array}\!\!\right)=\frac{t^{a_2}(1-t)^{\psi_p-1}_{+}}{\Gamma(\psi_p)}
{_2F_{1}}\!\left(\begin{array}{l}b_1-a_1,b_2-a_1
\\\psi_p\end{array}\!\!\vline\,1-t\right),
\end{equation}
where $a_1$ and $a_2$ may be interchanged on the right hand side.  If $\psi_p=-m$, $m=0,1,\ldots$, an easy calculation based on (\ref{eq:Meijer2reduced}) leads to
$$
G^{2,0}_{2,2}\!\left(\!t~\vline\begin{array}{l}b_1,b_2
\\a_1,a_2\end{array}\!\!\right)\!=\!\frac{t^{a_2}(b_1-a_1)_{m+1}(b_2-a_1)_{m+1}}{(m+1)!}
{_2F_{1}}\!\left(\begin{array}{l}b_1-a_1+m+1,b_2-a_1+m+1
\\m+2\end{array}\!\!\vline\,1-t\right),
$$
where again $a_1$ and $a_2$ may be interchanged.  The last formula holds for $t\in(0,1)$, see Property~3.  Hence, in view of N{\o}rlund's formula for $g_n(\a_{[k]};\b)$ for $p=2$ (see \cite[page~12]{KPSIGMA}), identities  (\ref{eq:p+1Fp-negintpsi})-(\ref{eq:p-1Fp-negintpsi}) take the form:
\begin{multline*}
\frac{\Gamma(\a)}{\Gamma(\b)}{_{3}F_2}\!\left(\!\!\left.\!\!\begin{array}{c}\sigma,\a\\
\b\end{array}\right|-z\!\right)=\frac{\Gamma(a_2)}{(1+z)^{\sigma}}\sum\limits_{j=0}^{m}\frac{(b_1-a_1)_j(b_2-a_1)_j}{j!\Gamma(a_2+j-m)}
{_{2}F_1}\!\!\left(\left.\!\!\begin{array}{c}\sigma,j-m\\a_2+j-m\end{array}\right|\frac{z}{1+z}\!\right)
\\
+\frac{(b_1-a_1)_{m+1}(b_2-a_1)_{m+1}}{(m+1)!}\!\int\limits_{0}^{1}\frac{t^{a_2-1}}{(1+zt)^{\sigma}}
{_2F_{1}}\!\left(\begin{array}{l}b_1-a_1+m+1,b_2-a_1+m+1
\\m+2\end{array}\!\!\vline\,1-t\right)\!dt,
\end{multline*}
where $-\psi_2=a_1+a_2-b_1-b_2=m\in\N_0$.  Similarly,
\begin{multline*}
\frac{\Gamma(\a)}{\Gamma(\b)}{_{2}F_2}\!\left(\!\left.\!\!\begin{array}{c}\a\\
\b\end{array}\right|-z\!\right)=e^{-z}\Gamma(a_2)\sum\limits_{j=0}^{m}\frac{(b_1-a_1)_j(b_2-a_1)_j}{j!\Gamma(a_2+j-m)}
{_{1}F_1}\!\!\left(\left.\!\!\begin{array}{c}j-m\\
a_2+j-m\end{array}\right|z\!\right)
\\
+\frac{(b_1-a_1)_{m+1}(b_2-a_1)_{m+1}}{(m+1)!}
\!\int\limits_{0}^{1}e^{-zt}t^{a_2-1}{_2F_{1}}\!\left(\!\begin{array}{l}b_1-a_1+m+1,b_2-a_1+m+1\\m+2\end{array}\!\!\vline\,1-t\right)\!dt,
\end{multline*}
and, for $\psi_2=b_1+b_2-a=-m+1/2$,
\begin{multline*}
\frac{\sqrt{\pi}\Gamma(a)}{\Gamma(\b)}{_{1}F_2}\!\left(\!\left.\!\!\begin{array}{c}a\\
\b\end{array}\right|-z\!\right)=\sum\limits_{j=0}^{m}(-1)^j\frac{(1/2)_j(b_1-a)_{m-j}(b_2-a)_{m-j}}{(m-j)!}
{_{0}F_1}\!\!\left(\left.\!\!\begin{array}{c}-\\
1/2-j\end{array}\right|z\!\right)
\\
+\frac{(b_1-a)_{m+1}(b_2-a)_{m+1}}{(m+1)!}
\!\int\limits_{0}^{1}\cos\left(2\sqrt{zt}\right)t^{-1/2}{_2F_{1}}\!\left(\!\begin{array}{l}b_1-a+m+1,b_2-a+m+1\\m+2\end{array}\!\!\vline\,1-t\right)\!dt.
\end{multline*}
The first two formulas still hold with $a_1$ and $a_2$ interchanged. These representations are presumably new.

\section{Applications of the integral representations}

\subsection{Inverse factorial series for ${}_{p+1}F_{p}$}

By factoring the generalized Stieltjes transform  (\ref{eq:Frepr}) into repeated Laplace transforms according to \cite[Theorem~8]{KPJCA} and applying (\ref{eq:qFqLaplace}),
we obtain (see also \cite[Theorem~4]{KarpJMS}):
\begin{multline*}
\frac{1}{z^m}{}_{p+1}F_{p}(m,\a;\b;-1/z)=\frac{\Gamma(\b)}{\Gamma(m)\Gamma(\a)}\int_{0}^{\infty}e^{-zu}u^{m-1}du\!\!\int\nolimits_{0}^{1}\!e^{-ux}G^{p,0}_{p,p}\left(x\left|\begin{array}{l}\!\!\b\!\!\\\!\!\a\!\!\end{array}\right.\right)\frac{dx}{x}
\\
=\frac{1}{\Gamma(m)}\!\int_{0}^{\infty}\!e^{-zu}u^{m-1}{_pF_p}(\a;\b;-u)du
=\frac{1}{\Gamma(m)}\int_0^1t^{z-1}(-\log{t})^{m-1}{}_{p}F_{p}(\a;\b;\log{t})dt
\\
=\frac{(-1)^{m-1}}{\Gamma(m)}\int_0^1t^{z-1}
\left(\sum\limits_{j=0}^{\infty}\frac{(\a)_j}{(\b)_jj!}(\log{t})^{m-1+j}\right)dt.
\end{multline*}
This formula is valid for any $m>0$, although for our purposes we only need to confine ourselves to $m\in\N$.
Further, according to \cite[Theorem~8.3]{Charalambides} we have
$$
(\log{t})^{m-1+j}=(m-1+j)!\sum\limits_{n=m-1+j}^{\infty}s(n,m-1+j)\frac{(t-1)^{n}}{n!},
$$
where $s(n,k)$ stands for the Stirling number of the first kind \cite[Section~8.2]{Charalambides}. Substituting this into the integrand above, we get:
\begin{multline*}
\sum\limits_{j=0}^{\infty}\frac{(\a)_j}{(\b)_jj!}(\log{t})^{m-1+j}
=\sum\limits_{j=0}^{\infty}\frac{(\a)_j(m-1+j)!}{(\b)_jj!}\sum\limits_{n=m-1+j}^{\infty}s(n,m-1+j)\frac{(t-1)^{n}}{n!}
\\
=\sum\limits_{n=m-1}^{\infty}\frac{(t-1)^{n}}{n!}\sum\limits_{j=0}^{n-m+1}\frac{(\a)_j(m-1+j)!}{(\b)_jj!}s(n,m-1+j)
=\sum\limits_{n=m-1}^{\infty}b_n(1-t)^{n},
\end{multline*}
where
$$
b_n:=\frac{(-1)^n}{n!}\sum\limits_{j=0}^{n-m+1}\frac{(\a)_j(m-1+j)!}{(\b)_jj!}s(n,m-1+j),
$$
and the series converges in the disk $|1-t|<1$.  Convergence follows from the fact that the repeated series on the right hand side of the first equality is easily seen to be absolutely convergent for $|1-t|<1$.  Substitution yields:
\begin{multline*}
\frac{1}{z^m}{}_{p+1}F_{p}(m,\a;\b;-1/z)=\frac{(-1)^{m-1}}{\Gamma(m)}\int_0^1t^{z-1}\left[\sum\limits_{n=m-1}^{\infty}b_n(1-t)^{n}\right]dt
\\
=\frac{(-1)^{m-1}}{\Gamma(m)}\sum\limits_{n=m-1}^{\infty}b_n\int_0^1t^{z-1}(1-t)^{n}dt=\frac{(-1)^{m-1}}{\Gamma(m)}\sum\limits_{n=m-1}^{\infty}\frac{b_nn!}{(z)_{n+1}}.
\end{multline*}
The inverse factorial series on the right converges for $\Re(z)>0$. This follows from the absolute convergence of the integral
on the right hand side of the first equality or from the general theory of inverse factorial series, see \cite[Theorems~III and IV]{Norlund14} and \cite[\S94~I,II]{Nielsen}.  The idea of inverse factorial series expansion of Stieltjes transform is also contained in the survey \cite{Weniger2010} by Weniger.  Rewriting the above formula with $w=1/z$ we arrive at the following theorem.

\begin{theorem}\label{th:invfact}
For arbitrary complex vectors $\a$, $\b$ and $m\in\N$, the inverse factorial expansion
$$
\frac{{}_{p+1}F_{p}(m,\a;\b;-w)}{\Gamma(\b)}=\frac{(-1)^{m-1}}{\Gamma(m)w^m}
\sum\limits_{n=m-1}^{\infty}\frac{(-1)^n}{(1/w)_{n+1}}\sum\limits_{j=0}^{n-m+1}\frac{(\a)_j(m-1+j)!}{\Gamma(\b+j)j!}s(n,m-1+j)
$$
converges for $\Re{w}>0$.
\end{theorem}

\textbf{Remark.} For general ${}_{p+1}F_{p}(\a';\b';-w)$ we can use the above theorem with $m=1$ by writing $\a=(1,\a')$, $\b=(1,\b')$.  We prefer to formulate it for general natural $m$ as the presence of a positive integer among the components of $\a'$ eliminates the need to extend the vectors $\a'$ and $\b'$.

\textbf{Remark.} In the recent preprint \cite{Costin} O.\:Costin and R.D.\:Costin introduced an extension of inverse factorial series convergent in domains larger than half-planes.

\subsection{Zeros of the Kummer and Bessel type functions}

In this section we show how (\ref{eq:qFqLaplace})  and (\ref{eq:q-1Fqcosine}) can be used to draw certain conclusions about zeros of ${}_{p}F_{p}$ and ${}_{p-1}F_{p}$ and derive some bounds for the latter for negative argument.  We start with an auxiliary fact that might be of independent interest.

\begin{lemma} \label{lm:Gdecrease}
Suppose that ${a_k}\le\min\{0,b_s-1\}$ for some indexes $k,s\in\{1,\ldots,p\}$ and $v_{\a_{[k]},\b_{[s]}}(t)$ defined in \emph{(\ref{eq:vge0})} is nonnegative on $[0,1]$ \emph{(}this holds, in particular, if $\b_{[s]}+\alpha\prec^W\a_{[k]}+\alpha$ for some $\alpha\in\R$\emph{)}. Then the function
$$
t\to G^{p,0}_{p,p}\left(t\left|\begin{array}{l}\!\!\b\\\!\!\a\end{array}\!\!\right.\right)
$$
is positive and decreasing on $(0,1)$.
\end{lemma}
\textbf{Proof}.  Set $\gamma=a_k\le0$, $\beta=b_s-a_k\ge1$, $\eta=-\beta-\gamma$,  $\b'=\b_{[s]}$, $\a'=\a_{[k]}$.  Then according to \cite[2.24.2.2]{PBM3} and in view of Properties~2 and 3 found in the Appendix to this paper, we have:
\begin{multline*}
G^{p,0}_{p,p}\left(x\left|\begin{array}{l}\!\b\\\!\a\end{array}\!\!\right.\right)
=G^{p,0}_{p,p}\left(x\left|\begin{array}{l}\!\!\beta+\gamma,\b'+\beta+\eta+\gamma\\\gamma, \a'+\beta+\eta+\gamma\end{array}\!\!\right.\right)
\\
=x^{\gamma}G^{p,0}_{p,p}\left(x\left|\begin{array}{l}\!\!\beta,\b'+\beta+\eta\\0, \a'+\beta+\eta\end{array}\!\!\right.\right)
=
\frac{x^{\gamma}}{\Gamma(\beta)}\int_{x}^1t^{\eta}(t-x)^{\beta-1}G^{p-1,0}_{p-1,p-1}\left(t\left|\begin{array}{l}\!\!\b' \\ \!\!\a'\end{array}\!\!\right.\right)dt.
\end{multline*}
By the hypotheses of the lemma $v_{\a',\b'}(t)\ge0$ on $[0,1]$, so that by Property~9 (see Appendix), the $G$ function in the integrand is nonnegative.  Combined with the conditions $\beta\ge1$ and $\gamma\le0$ this implies that the rightmost term in the above chain is decreasing and so does the leftmost term.  $\hfill\square$

\begin{theorem}\label{th:pFpnozeros}
Let $\a$, $\b$ be positive vectors. Suppose that $a_k\le\min\{1,b_s-1\}$ for some indexes $k,s\in\{1,\ldots,p\}$ and $v_{\a_{[k]},\b_{[s]}}(t)\ge0$ on $[0,1]$ \emph{(}in particular, $\b_{[s]}\prec^W\a_{[k]}$  is sufficient\emph{)}.  Then ${}_{p}F_p(\a;\b;z)$ has no real zeros and all its zeros lie in the open right half plane $\Re(z)>0$.
\end{theorem}
\textbf{Proof}. Indeed, under the hypotheses of the theorem formula (\ref{eq:qFqLaplace}) is applicable. After the change of variable $t=1-u$ we get:
$$
e^{z}{}_{p}F_p\left.\left(\begin{matrix}\a\\\b\end{matrix}\right\vert -z\right)=
\frac{\Gamma(\b)}{\Gamma(\a)}\int_0^1e^{zu}G^{p,0}_{p,p}\left(1-u\left|\begin{array}{l}\!\!\b-1\!\!\\\!\!\a-1\!\!\end{array}\right.\right)\!du.
$$
By Lemma~\ref{lm:Gdecrease}, the $G$ function in the integrand is positive and increasing on $(0,1)$ and is clearly not a step function. The claim now follows by \cite[Theorem~2.1.7]{Sedlet2005}. $\hfill\square$

\begin{theorem}\label{th:p-1Fpestimate}
Let $\a'\in\R^{p-1}$, $\b\in\R^{p}$ be positive vectors. Set $\a=(\a',1/2)$ and assume that $v_{\a,\b}(t)\ge0$ on $[0,1]$ \emph{(}it suffices that $\b\prec^W\a$\emph{)}.
Then for any $x>0$,
$$
\left|{}_{p-1}F_p\left.\left(\begin{matrix}\a'\\\b\end{matrix}\:\right\vert -x\right)\right|<1.
$$
In particular, the functions ${}_{p-1}F_p\left(\a';\b; x\right)\pm1$ have no real zeros other than $x=0$ \emph{(}in case of the ''minus''  sign\emph{)}.
\end{theorem}
\textbf{Proof.}  Set $t=u^2$ and replace $z\to{z^2/4}$ in (\ref{eq:q-1Fqcosine}) to get:
\begin{equation}\label{eq:p-1Fpcos}
{}_{p-1}F_p\left.\left(\begin{matrix}\a'\\\b\end{matrix}\:\right\vert -z^2/4\right)
=\frac{2\Gamma(\b)}{\sqrt{\pi}\Gamma(\a')}\int_0^1\cos(zu)G^{p,0}_{p,p}\left(u^2\left|\begin{array}{l}\!\!\b-1/2\!\!\\\!\!\a'-1/2,0\!\!\end{array}\right.\right)\!du.
\end{equation}
From this formula, for any real $z$ we obtain the estimate (recall that $\a=(\a',1/2)$)
$$
\left|{}_{p-1}F_p\left.\left(\begin{matrix}\a'\\\b\end{matrix}\:\right\vert -z^2/4\right)\right|\leq\frac{2\Gamma(\b)}{\Gamma(\a)}\int_0^1\left|G^{p,0}_{p,p}\left(u^2\left|\begin{array}{l}\!\!\b-1/2\!\!\\\!\!\a-1/2\!\!\end{array}\right.\right)\right|\!du=1,
$$
where the last equality follows from nonnegativity of the $G$ function as indicated in Property~9 in the Appendix of this paper.  Hence, to compute the last integral we can drop the absolute value, substitute $t=u^2$ and  use (\ref{eq:GMellin}).  The inequality is in fact strict for all real $z\ne0$ as can be seen from (\ref{eq:p-1Fpcos}) by the mean value theorem.  $\hfill\square$

\begin{theorem}\label{th:p-1Fpsine}
Let $\a'\in\R^{p-1}$, $\b\in\R^{p}$ be positive vectors.  Set $\a=(\a',3/2)$ and suppose that $a_k\le\min\{1,b_s-1\}$ for some indexes $k,s\in\{1,\ldots,p\}$ while $v_{\a_{[k]},\b_{[s]}}(t)\ge0$ on $[0,1]$ \emph{(}or $\b_{[s]}\prec^W\a_{[k]}$\emph{)}.  Then $0<{}_{p-1}F_p(\a';\b;x)<1$ for all $x<0$.  In particular,
${}_{p-1}F_p(\a';\b;x)$ has no real zeros.
\end{theorem}
\textbf{Proof.} Note first that ${}_{p-1}F_p(\a';\b;x)$ has no zeros for $x\ge0$ as it is obvious from the series  representation (\ref{eq:pFqdefined}).  Second, a short reflection shows that the hypotheses of this theorem imply the hypotheses of Theorem~\ref{th:p-1Fpestimate} and hence ${}_{p-1}F_p(\a';\b;x)<1$. It remains to prove that  ${}_{p-1}F_p(\a';\b;x)$ has no zeros for negative $x$.  This fact will be derived from the following representation:
\begin{equation}\label{eq:p-1Fpsine}
z{}_{p-1}F_p\left.\left(\begin{matrix}\a'\\\b\end{matrix}\:\right\vert -z^2/4\right)
=\frac{4\Gamma(\b)}{\sqrt{\pi}\Gamma(\a')}\int_0^1\sin(zu)
G^{p,0}_{p,p}\left(u^2\left|\begin{array}{l}\!\!\b-1\!\!\\\!\!\a'-1,1/2\!\!\end{array}\right.\right)\!du.
\end{equation}
To prove (\ref{eq:p-1Fpsine}) substitute $t=u^2$ and exchange the order of summation and integration to get:
\begin{multline*}
\int_0^1\sin(zu)G^{p,0}_{p,p}\left(u^2\left|\begin{array}{l}\!\!\b-1\!\!\\\!\!\a'-1,1/2\!\!\end{array}\right.\right)\!du
=\int_0^1\sin(z\sqrt{t})G^{p,0}_{p,p}\left(t\left|\begin{array}{l}\!\!\b-3/2\!\!\\\!\!\a'-3/2,0\!\!\end{array}\right.\right)\!\frac{dt}{2}
\\
=\sum\limits_{k=0}^{\infty}\frac{(-1)^kz^{2k+1}}{2(2k+1)!}\int_0^1t^{k-1}G^{p,0}_{p,p}\left(t\left|\begin{array}{l}\!\!\b\!\!\\\!\!\a',3/2\!\!\end{array}\right.\right)\!dt
=\frac{z}{2}\sum\limits_{k=0}^{\infty}\frac{(-z^2/4)^{k}}{(3/2)_kk!}\frac{\Gamma(\a'+k)\Gamma(3/2+k)}{\Gamma(\b+k)}
\\
=\frac{z\Gamma(\a')\Gamma(3/2)}{2\Gamma(\b)}\sum\limits_{k=0}^{\infty}\frac{(\a')_k(-z^2/4)^{k}}{(\b)_kk!}=\frac{z\sqrt{\pi}\Gamma(\a')}{4\Gamma(\b)}
{}_{p-1}F_p\left.\left(\begin{matrix}\a'\\\b\end{matrix}\:\right\vert -z^2/4\right).
\end{multline*}
According to Lemma~\ref{lm:Gdecrease} (applied with $\a$, $\b$ replaced with $\a-1$ and $\b-1$), the hypotheses of the theorem imply that the $G$ function in the integrand of (\ref{eq:p-1Fpsine}) is positive and decreasing.  We are now in the position to apply \cite[Theorem~2.1.5]{Sedlet2005} (the English translation has am important omission in the formulation, so we prefer to refer to the Russian original), which states that $\int_0^1\sin(zu)f(u)du$ has no real zeros apart from $z=0$ if $f$ is positive and decreasing on (0,1) and is not a step function with rational jump points. Clearly, the $G$ function is not a step function on $(0,1)$ since it is a combination of powers and analytic functions by (\ref{eq:Gp0qp}) and the claim follows.$\hfill\square$

\textbf{Remark}. By renumbering parameters any function satisfying Theorem~\ref{th:p-1Fpsine} can be written as
${}_{p-1}F_{p}\!\left(\alpha,\a;\beta_1,\beta_2,\b;x\right)$, where  $0<\alpha\le1$, $\beta_1\ge\alpha+1$,  $\beta_2\ge3/2$,
$\a>0$ and $v_{\a,\b}(t)\ge0$ on $[0,1]$.

\textbf{Remark}.  For arbitrary $j\in\{1,\ldots,p\}$ we can use the representation \cite[Remark~on~page~124]{KarpJMS}
\begin{multline*}
{_{p-1}F_p}(\a;\b;-z)=\frac{\Gamma(\b_{[j]})}{\Gamma(\a)}
\int\limits_{0}^{1}{_0F_{1}}(-;b_j;-zt)
G^{p-1,0}_{p-1,p-1}\left(t\left|\begin{array}{l}\!\!\b_{[j]}\!\!\\\!\!\a\!\!\end{array}\right.\right)\frac{dt}{t}
\\
=\frac{z^{(1-b_j)/2}\Gamma(\b)}{\Gamma(\a)}
\int\limits_{0}^{1}J_{b_j-1}(2\sqrt{zt})
G^{p-1,0}_{p-1,p-1}\left(t\left|\begin{array}{l}\!\!\b_{[j]}-(b_j+1)/2\!\!\\\!\!\a-(b_j+1)/2\!\!\end{array}\right.\right)dt
\end{multline*}
to improve Theorem~\ref{th:p-1Fpestimate}. Here $J_{\nu}$ is the Bessel function of the first kind. Using $|J_{\nu}(x)|\leq1$ for all real $x$ if $\nu\ge0$ \cite[10.14.1]{NIST} we obtain:
\begin{equation}\label{p-1Fpbesselbound}
\left|{_{p-1}F_p}(\a;\b;-x)\right|\leq\frac{x^{-b_j/2+1/2}\Gamma(\b)\Gamma(\a-(b_j-1)/2)}{\Gamma(\b_{[j]}-(b_j-1)/2)\Gamma(\a)}~\text{for}~x>0,
\end{equation}
if $b_j\geq1$, $\a-(b_j+1)/2>0$ and $\b_{[j]}\prec^W\a$ for some $j\in{1,\ldots,p}$.  The constant can be further improved
by employing the result of Landau \cite{Landau}: $|J_{\nu}(x)|\leq\alpha\nu^{-1/3}$ for all real $x$ if $\nu\ge0$, where $\alpha\approx0.674885$.  This gives an improvement over (\ref{p-1Fpbesselbound}) if $b_j\ge1.31$.  It follows from Properties 5 and 6 (see Appendix) that we can relax the conditions on parameters to $\a-(b_j+1)/2>0$ and $\sum_{k\ne{j}}b_k-\sum_{k}a_k>0$ at the price of losing the exact expression for the constant in (\ref{p-1Fpbesselbound}), i.e. we get a bound of the form
$Cx^{(1-b_j)/2}$. Furthermore, we can use another bound due to Landau  in \cite{Landau}: $|J_{\nu}(x)|\leq\beta|x|^{-1/3}$
valid for all real $x$ and $\nu\ge0$ with $\beta\approx0.785747$.  This allows the reduction of the power factor in (\ref{p-1Fpbesselbound}) to $x^{-b_j/2+1/3}$, at the expense of slightly increasing the constant factor.  These bounds can also be combined with the well known estimate $|J_{\nu}(x)|\leq2^{-\nu}x^{\nu}/\Gamma(\nu+1)$ (see \cite[10.14.4]{NIST}), valid for $\nu\ge-1/2$, to get improved inequalities for $|{}_{p-1}F_p(\a;\b;-x)|$ in different $x$ regions.

\textbf{Remark.} For positive $x$, a two-sided bound for ${_{p-1}F_p}(x)$ was found in \cite[Theorems~10, 11]{KarpJMS}.

\subsection{Radial positive definite functions}

The purpose of this section is to demonstrate that the generalized hypergeometric functions provide a plethora of examples of radial positive definite functions
well suited for formulating and/or verifying hypotheses about such functions. It is worth mentioning that hypergeometric examples of radial positive definite functions have been considered recently in \cite{PorcuZast}. Let us remind the reader that a continuous function $f$ on $(0,\infty)$ is called $n$-RPDF (radial positive definite in dimension $n$) if for each $m\in\N$
$$
\sum\limits_{i,j=1}^{m}f(\|t_i-t_j\|_n)\xi_i\overline{\xi}_j\ge0,~~\forall\{t_1,\ldots,t_m\}\subset\R^n,
~~\forall\{\xi_1,\ldots,\xi_m\}\subset\C.
$$
The class of $n$-RPD functions is denoted by $\Phi_n$. The above definition and many further details can be found, for instance, in the two recent papers \cite{GMO1,GMO2}.
The class $\Phi_n$ has been characterized by Schoenberg in 1938:  $f\in\Phi_n$ with $f(0)=1$ iff $f(r)=\int_0^\infty \Omega_n(rt)\nu_f(dt)$,
where $\nu_f$ is a probability measure on $[0,\infty)$ uniquely determined by $f$ and
$$
\Omega_n(s):=\sum\limits_{j=0}^{\infty}\frac{(-s^2/4)^j}{(n/2)_jj!}={}_0F_{1}(-;n/2;-s^2/4).
$$
Classes $\Phi_n$ are known to be nested: $\Phi_{n+1}\subset\Phi_n$, and the inclusion is proper. The class
$$
\Phi_{\infty}:=\bigcap\limits_{n\ge1}\Phi_n
$$
has been also characterized by Schoenberg as follows: $f\in\Phi_{\infty}$ with $f(0)=1$ iff $f(r)=\int_0^\infty e^{-tr^2}\nu_f(dt)$, where $\nu_f$ is a probability measure on $[0,\infty)$.  These characterizations allow us to give sufficient conditions on the parameters of the generalized hypergeometric functions which guarantee that they are radial positive definite for a certain dimension.  The results for the Gauss and Kummer type functions are simple and complete.
\begin{theorem}\label{th:GaussKummerRPDF}
Suppose $\sigma,\a,\b>0$ and $v_{\a,\b}(t)\ge0$ on $[0,1]$. Then
$$
{_{p}F_{p}}\left(\left.\!\!\begin{array}{c}\a\\\b\end{array}\right|-r^2\!\right)\in\Phi_{\infty}~~\text{and}~~{_{p+1}F_{p}}\left(\left.\!\!\begin{array}{c}\sigma,\a\\\b\end{array}\right|-r^2\!\right)\in\Phi_{\infty}.
$$
\end{theorem}
\textbf{Proof.} Indeed, the representation
$$
{_{p}F_{p}}\left(\left.\!\!\begin{array}{c}\a\\\b\end{array}\right|-r^2\!\right)=\frac{\Gamma(\b)}{\Gamma(\a)}\int_0^1e^{-tr^2}G^{p,0}_{p,p}\!\left(\!t~\vline\begin{array}{l}\b\\\a\end{array}\!\!\right)\frac{dt}{t}
$$
is a rewriting of (\ref{eq:qFqLaplace}), while the representation
$$
{_{p+1}F_{p}}\left(\left.\!\!\begin{array}{c}\sigma,\a\\\b\end{array}\right|-r^2\!\right)=\frac{\Gamma(\b)}{\Gamma(\a)\Gamma(\sigma)}\int_0^{\infty}e^{-tr^2}G^{p+1,0}_{p,p+1}\!\left(\!t~\vline\begin{array}{l}\b\\\sigma,\a\end{array}\!\!\right)\frac{dt}{t}
$$
is given in \cite[formula (10)]{KarpJMS}. The function $G^{p,0}_{p,p}$ in the first representation is nonnegative by Property~9.  Nonnegativity of the weight function in the second representation follows from the formula \cite[2.24.3.1]{PBM3}:
$$
G^{p+1,0}_{p,p+1}\left(t\left|\begin{array}{l}\!\!\b\!\!\\\!\!\sigma,\a\!\!\end{array}\right.\right)=t^{\sigma}
\int\limits_{1}^{\infty}e^{-ty}y^{\sigma-1}G^{p,0}_{p,p}\left(\frac{1}{y}\left|\begin{array}{l}\!\!\b\!\!\\\!\!\a\!\!\end{array}\right.\right)dy.\hspace{2cm}\square
$$
For the Bessel type functions we begin with the following monotonicity theorem.
\begin{theorem}
Suppose that $0<\a'\le\a$, $\b'\ge\b>0$ \emph{(}understood element-wise\emph{)} and ${_{p-1}F_{p}}\left(\a;\b;-r^2\right)\in\Phi_{n}$ for some $n\in\N$.  Then
${_{p-1}F_{p}}\left(\a';\b';-r^2\right)\in\Phi_{n}$.
\end{theorem}
\textbf{Proof.} As linear combination of positive definite matrices is positive definite, it is easy to see from the definition of RPDF that  $g(y)=\int_0^1f(yx)k(x)dx\in\Phi_{n}$ if $f\in\Phi_n$ and $k(x)\ge0$.
Assume that $\a'=(a_1',a_2,\ldots,a_p)$ with $0<a_1'<a_1$ and $\b'=\b$. Then an easy calculation using termwise integration yields:
$$
{_{p-1}F_{p}}\left(\a';\b';-r^2\right)=\frac{2}{B(a_1',a_1-a_1')}\int_0^1t^{2a_1'-1}(1-t^2)^{a_1-a_1'-1}{_{p-1}F_{p}}\left(\a;\b;-(rt)^2\right)dt,
$$
and the claim follows. General $\a'\le\a$ and $\b'\ge\b$ can be treated similarly, taking each pair of non equal components one by one. $\hfill\square$

\begin{theorem}\label{th:BesselRPDF}
Suppose that $\a,\b>0$. If ${_{p-1}F_{p}}\left(\a;\b;-r^2\right)\in\Phi_{n}$, then
$$
\psi_p=\sum\limits_{j=1}^{p}b_j-\sum\limits_{j=1}^{p-1}a_j\ge\frac{n}{2}.
$$
The Schoenberg measure is supported on $[0,2]$ and is given by
$$
\nu(dt)=\frac{2\Gamma(\b)}{\Gamma(\a)\Gamma(n/2)}G^{p,0}_{p,p}\!\left(\frac{t^2}{4}~\vline\begin{array}{l}\b\\\a,n/2\end{array}\!\!\right)\frac{dt}{t}
$$
for $\psi_p>n/2$ and by $\widetilde{\nu}(dt)=\nu(dt)+\delta_2$, where $\delta_2$ is the Dirac measure concentrated at $2$, for $\psi_p=n/2$. In particular, ${_{p-1}F_{p}}\left(\a;\b;-r^2\right)\notin\Phi_{\infty}$ for any $\a,\b>0$.

Conversely, set $\a'=(n/2,\a)$ and suppose $v_{\a',\b}(t)\ge0$ on $[0,1]$ \emph{(}which implies $\psi_p\ge{n/2}$\emph{)}. Then ${_{p-1}F_{p}}\left(\a;\b;-r^2\right)\in\Phi_{n}$.
If $\psi_p={n/2}$, then ${_{p-1}F_{p}}\left(\a;\b;-r^2\right)\in\Phi_n\setminus\Phi_{n+1}$.
\end{theorem}
\textbf{Proof.} Indeed, suppose that ${_{p-1}F_{p}}\left(\a;\b;-r^2\right)=\int_0^{\infty}\Omega_n(rt)\nu(dt)$. The argument used in the proof of \cite[Theorem~2.1.1]{LinnikOstov}
shows \emph{mutatis mutandis} that existence of all derivatives at $r=0$ of the function on the left hand side implies that all moments of the measure $\nu(dt)$ are finite.
Power series (\ref{eq:pFqdefined}) and definition of $\Omega_n$ then show that the even moments of $\nu$ are given by
$$
\nu_{2k}=\int\limits_{0}^{\infty}t^{2k}\nu(dt)=4^k\frac{(\a)_k(n/2)_k}{(\b)_k},~~~k=0,1,\ldots
$$
Next, changing variable $t=2u$ in the above integral, we get
$$
\nu_{2k}=\int\limits_{0}^{\infty}(2u)^{2k}\widetilde{\nu}(du)=4^k\frac{(\a)_k(n/2)_k}{(\b)_k}~~\text{or}~~\widetilde{\nu}_{2k}=\int\limits_{0}^{\infty}u^{2k}\widetilde{\nu}(du)=\frac{(\a)_k(n/2)_k}{(\b)_k},
$$
where $\widetilde{\nu}$ is the image measure of $\nu$ under this change of variable.  If $\psi_p<n/2$, then the sequence $\widetilde{\nu}_{2k}\sim k^{n/2-\psi_p}$ as $k\to\infty$.  However, it is easy to see that a Stieltjes moment sequence, either tends to the atom of the representing measure at 1 (if the support of the measure is contained in $[0,1]$) or grows at least exponentially (if the support of the measure contains points outside $[0,1]$), see a related result in \cite[Lemma~2.9]{BergDuran}.  Therefore, the sequence $(\a)_k(n/2)_k/(\b)_k$ is not a Stieltjes moment sequence for $\psi_p<n/2$, proving our first claim.  For $\psi_p>n/2$ the expression for $\nu(dt)$ follows from formula (\ref{eq:GMellin}) found in Property~7 contained in the Appendix.  If $\psi_p=n/2$ the expression for $\widetilde{\nu}(dt)$
follows from Theorem~\ref{th:zerobalanced} or from (\ref{eq:GMellinNorlund}).

The first claim in the converse statement for $\psi_p>n/2$ is immediate from the representation (verified by termwise integration)
$$
{_{p-1}F_{p}}\left(\left.\!\!\begin{array}{c}\a\\\b\end{array}\right|-r^2\!\right)
=\frac{2\Gamma(\b)}{\Gamma(\a)\Gamma(n/2)}\int_0^2\Omega_n(rt)G^{p,0}_{p,p}\!\left(\frac{t^2}{4}~\vline\begin{array}{l}\b\\\a,n/2\end{array}\!\!\right)\frac{dt}{t}
$$
and Property~9 given in the Appendix. To establish the second claim,  ${_{p-1}F_{p}}\left(\a;\b;-r^2\right)\in\Phi_n\setminus\Phi_{n+1}$ for $\psi_p=n/2$, we invoke \cite[Theorem~3.1]{GMO2} which implies that a function belongs to $\Phi_n\setminus\Phi_{n+1}$ if its Schoenberg measure contains an atom. By Theorem~\ref{th:zerobalanced} the representing measure of ${_{p-1}F_{p}}$ indeed has an atom at $t=2$ when $\psi_p=n/2$. $\hfill\square$

\textbf{Remark.} Leonid Golinskii observed that the membership $[\Omega_n(r)]^2\in\Phi_{2n-1}\!\setminus\!\Phi_{2n}$ proved in \cite[Theorem~1.3(ii)]{GMO2} also follows directly from Theorem~\ref{th:BesselRPDF}. Indeed, using the formula for the product of Bessel functions \cite[Exercise~16, p.237]{AAR} we conclude that
$$
[\Omega_n(r)]^2={}_2F_{3}\left(\left.\!\!\begin{array}{c}n/2,n/2-1\\n/2,n/2,n-1\end{array}\right|-r^2\!\right).
$$
Here $\psi_3=(2n-1)/2$ so that by the last statement of Theorem~\ref{th:BesselRPDF} $[\Omega_n(r)]^2\in\Phi_{2n-1}\!\setminus\!\Phi_{2n}$.

\textbf{Remark.} We mention further connections of hypergeometric functions with two more functional classes considered in \cite{GMO1}.  For nonnegative, monotone decreasing functions $f$, normalized by $f(0)=1$ the authors proved that the Schoenberg operator associated with the matrix $[f(\|t_i-t_j\|_n)]_{i,j=1}^{m}$, $m\le\infty$, is bounded on $l^2$ under the additional condition $t^{d-1}f\in{L_1(\R_+)}$, where $d$ is the dimension of the linear span of $t_1,\ldots,t_m$.  Representations for the Gauss and Kummer type functions
${}_{p+1}F_p(\a;\b;-r^2)$ and ${}_{p}F_p(\a;\b;-r^2)$ exhibited in the proof of Theorem~\ref{th:GaussKummerRPDF} show that both are nonnegative, monotone decreasing and properly normalized.  Asymptotic formulas \cite[16.11.6,16.11.7]{NIST} show that the condition $t^{d-1}f\in{L_1(\R_+)}$ is satisfied for dimensions $d<\min(\a)$.  Finally, \cite[Theorem~1.7]{GMO1} shows that under the conditions of Theorem~\ref{th:GaussKummerRPDF}, the functions ${}_{p+1}F_p(\a;\b;-r^2)$ and ${}_{p}F_p(\a;\b;-r^2)$  are \emph{strongly} $X$ positive definite for any separated set $X$.  The definition of strongly positive definite functions is found in \cite[Definition~1.5]{GMO1}, their importance is also explained in \cite{GMO1} and references therein.

\section{General complex parameters}

\subsection{Representations of GHF}

The following decomposition is straightforward:
\begin{equation}\label{eq:pFq-decompose}
{}_{p}F_q\left.\left(\begin{matrix}\a\\\b\end{matrix}\right\vert z\right)=
\sum_{k=0}^{n-1}\frac{(\a)_k}{(\b)_kk!}z^k
+\frac{(\a)_nz^{n}}{(\b)_nn!}{}_{p+1}F_{q+1}\left.\left(\begin{matrix}\a+n,1\\\b+n,n+1\end{matrix}\right\vert z\right).
\end{equation}
The new parameter vectors $\a'=(\a+n,1)$ and $\b'=(\b+n,n+1)$ clearly satisfy $\Re(\a')>0$ and for  $p=q$ also $\Re(\psi_p')>0$ for sufficiently large $n$, where $\psi_p':=\sum_{k=1}^{p+1}(b_k'-a_k')$.  This observation immediately leads to

\begin{theorem}\label{th:pFqdecompose}
For arbitrary $\a,\b\in\C^{p}$ choose an $n\in\N_0$ satisfying $\Re(\a)+n>0$ and $\Re(\psi_p)+n>0$, where $\psi_p:=\sum_{k=0}^{p}(b_k-a_k)$.  Then
\begin{equation}\label{eq:G-Stieltjes}
\frac{1}{\Gamma(\b)}{}_{p+1}F_p\left.\left(\begin{matrix}\sigma,\a\\\b\end{matrix}\right\vert-z\right)=
\sum_{k=0}^{n-1}\frac{(\sigma)_k(\a)_k(-z)^k}{\Gamma(\b+k)k!}
+{(\sigma)_n(-z)^n\over\Gamma(\a)}\int_0^1\frac{\widetilde{G}_n(t)dt}{(1+zt)^{\sigma+n}}
\end{equation}
for all $\sigma\in\C$ and $z\in\C\!\setminus\!(-\infty,-1]$, and
\begin{equation}\label{eq:G-Laplace}
\frac{1}{\Gamma(\b)}{}_{p}F_p\left.\left(\begin{matrix}\a\\\b\end{matrix}\right\vert-z\right)=
\sum_{k=0}^{n-1}\frac{(\a)_k(-z)^k}{\Gamma(\b+k)k!}+{(-z)^n\over\Gamma(\a)}
\int_0^1e^{-zt}\widetilde{G}_n(t)dt
\end{equation}
for all $z\in\C$, where
\begin{equation}\label{eq:tildeGn-defined}
\widetilde{G}_n(t):=G^{p+1,0}_{p+1,p+1}\left(t\left|\begin{array}{l}\!\!\b-1+n,n\!\!\\\!\!\a-1+n,0\!\!\end{array}\right.\right).
\end{equation}
If $\a,\b\in\R^{p}$, then there exists $n\in\N_0$ such that $\widetilde{G}_n(t)\ge0$ for $t\in(0,1)$.  In particular, this $n$ can be chosen from the condition $(\b+n,n+1)\prec^W(\a+n,1)$.

Furthermore, for $\a\in\C^{p-1}$, $\b\in\C^p$ and $n\in\N_0$ satisfying $\Re(\a)+n>0$ and $\Re(\widehat{\psi_p})+2n-1/2>0$, where $\widehat{\psi_p}:=\sum_{k=1}^{p}b_k-\sum_{k=1}^{p-1}a_k$, we have
\begin{equation}\label{eq:G-cosine}
\frac{1}{\Gamma(\b)}{}_{p-1}F_p\left.\left(\begin{matrix}\a\\\b\end{matrix}\right\vert -z\right)
=\sum_{k=0}^{n-1}\frac{(\a)_k(-z)^k}{\Gamma(\b+k)k!}+{(-z)^n\over\sqrt{\pi}\Gamma(\a)}
\int_0^1\cos(2\sqrt{zt})\widehat{G}_n(t)dt,
\end{equation}
for all $z\in\C$, where
\begin{equation}\label{eq:hatGn-defined}
\widehat{G}_n(t):=G^{p+1,0}_{p+1,p+1}\left(t\left|\begin{array}{l}\!\!\b-1+n,n\!\!\\\!\!\a-1+n,-1/2,0\!\!\end{array}\right.\right). \end{equation}
If $\a\in\R^{p-1}$, $\b\in\R^{p}$, then there exists $n\in\N_0$ such that $\widehat{G}_n(t)\ge0$ for $t\in(0,1)$.  In particular, this $n$ can be chosen from the condition $(\b+n,n+1)\prec^W(\a+n,1/2,1)$.
\end{theorem}

\textbf{Proof.} For complex parameters satisfying the conditions of the theorem, formulas (\ref{eq:G-Stieltjes}), (\ref{eq:G-Laplace}) and (\ref{eq:G-cosine}) follow from (\ref{eq:pFq-decompose}) combined with (\ref{eq:Frepr}), (\ref{eq:qFqLaplace}) and (\ref{eq:q-1Fqcosine}).  Suppose now that $\a\in\R^{p}$, $\b\in\R^{p}$ are arbitrary.
To prove that $(\b+n,n+1)\!\prec^W\!(\a+n,1)$ for some $n\in\N_0$, assume that  $\sum_{j=1}^{k}a_j>\sum_{j=1}^{k}b_j$ for some $k$. Then, clearly,
$1+\sum_{j=1}^{k-1}a_j+(k-1)n\leq\sum_{j=1}^{k}b_j+kn$ for a sufficiently large $n$. The sum on the left has the form shown since $1\le\min(a_1+n,\ldots,a_p+n)$ for sufficiently large $n$.  Similarly, for $\a\in\R^{p-1}$, $(\b+n,n+1)\prec^W(\a+n,1/2,1)$ for sufficiently large $n$.  Nonnegativity of $\widetilde{G}_n$ and $\widehat{G}_n$ now follow  by Property~9 in the Appendix.  $\hfill\square$

\textbf{Remark.} As before the condition $\b'=(\b+n,n+1)\!\prec^W\!(\a+n,1)=\a'$ in the above theorem can be replaced by the weaker condition
$$
v_{\a',\b'}(t)=\sum\nolimits_{k=1}^{p}(t^{a_k+n}-t^{b_k+n})+t-t^{n+1}\ge0~\text{for}~t\in[0,1].
$$

\subsection{Regularization of the integrals containing $G^{p,0}_{p,p}$}

Decomposition formulas (\ref{eq:G-Stieltjes}), (\ref{eq:G-Laplace}) and  (\ref{eq:G-cosine}) can be viewed as manifestations of a more general phenomenon.
Define $\CB^{\infty}[0,1]$ to be the class of functions on $[0,1]$ that have bounded derivatives of all orders.  If $\varphi\in\CB^{\infty}[0,1]$ then the integral
\begin{equation}\label{eq:Gphiintegral}
\int_0^1\!\!G_0(t)\varphi(t)dt,~~\text{where}~~
G_0(t)=G^{p,0}_{p,p}\!\left(t\left\vert\begin{matrix}\b-1\\\a-1\end{matrix}\right.\right),
\end{equation}
converges (i.e. exists as an improper integral) if the next two conditions are satisfied:
\begin{equation}\label{eq:a-psi-defined}
\Re(\a)>0~~\text{and}~~\Re{\psi_p}=\Re\!\left[\sum\nolimits_{k=1}^{p}(b_k-a_k)\right]>0.
\end{equation}
This is implied by Properties~5 and 6 found in Appendix. Furthermore, according to (\ref{eq:Norlund1}), it exists if $\psi_p=0,-1,-2,\ldots$
The purpose of this section is to define a regularization of the integral (\ref{eq:Gphiintegral}) valid for arbitrary complex parameters. Choosing $\varphi$ to be the generalized Stieltjes or the exponential kernel will naturally lead representations of generalized hypergeometric functions  equivalent to (\ref{eq:G-Stieltjes}) and  (\ref{eq:G-Laplace}) above.  Curiously enough, taking $\varphi$ to be equal to the cosine Fourier kernel leads to the representation of the Bessel type function that is different from (\ref{eq:G-cosine}).

To convert the set $\CB^{\infty}[0,1]$ into a test function space, we introduce the following definition of convergence in $\CB^{\infty}[0,1]$:
the sequence $\varphi_j$  converges to an element $\varphi\in\CB^{\infty}[0,1]$ if
$$
\max\limits_{x\in[0,1]}|\varphi^{(k)}_j(x)-\varphi^{(k)}(x)|\to0~\text{as}~j\to\infty
$$
for each nonnegative integer $k$. This space can be viewed as a space of restrictions of smooth periodic functions (say with period $2$) considered in \cite[Chapter~3, paragraph~2]{Beals} to the interval $[0,1]$. Then it follows from \cite[Theorem~2.1]{Beals} that this space is complete.
\medskip

\noindent\textbf{Definition~1.} For arbitrary complex $\a$ and $\b$, $-\b\notin\N_0$,  choose a nonnegative integer $n$ satisfying $\Re(\a)+n>0$ and $\Re(\psi_p)+n>0$. Define a regularization of the integral (\ref{eq:Gphiintegral}) as the distribution $\G_0$ acting on a test function $\varphi\in\CB^{\infty}[0,1]$ according to the formula
\begin{equation}\label{eq:actionG0}
\langle\G_0,\varphi\rangle=\sum_{k=0}^{n-1}\frac{(\a)_k}{(\b)_kk!}\varphi^{(k)}(0)
+\frac{\Gamma(\b)}{\Gamma(\a)}\int_0^1\widetilde{G}_n(t)\varphi^{(n)}(t)dt,
\end{equation}
where
\begin{equation}\label{eq:tildeGn-defiend}
\widetilde{G}_n(t):=G^{p+1,0}_{p+1,p+1}\left(t\left|\begin{array}{l}\!\!\b-1+n,n\!\!\\\!\!\a-1+n,0\!\!\end{array}\right.\right), \hskip 1cm n=0,1,\ldots
\end{equation}
Clearly, $\widetilde{G}_0=G_0$ as defined in (\ref{eq:Gphiintegral}). Furthermore, if  $n=0$  the finite sum in (\ref{eq:actionG0}) is understood to be empty, so that (\ref{eq:actionG0}) reduces to  a multiple of (\ref{eq:Gphiintegral}).  The asymptotic behavior  of $\widetilde{G}_n(t)$ (as $t\to0$ and $t\to1$), contained in Properties~5 and 6 in the Appendix, shows that the integral in (\ref{eq:actionG0}) exists (as a finite number) for all $\varphi\in\CB^{\infty}[0,1]$ under the conditions stated in Definition~1.  Note that if $\Re(\psi_p)>0$, the function $\widetilde{G}_n(t)$ can be computed as the $n-$th primitive of $G_0(x)$ that satisfies $\widetilde{G}_n^{(k)}(1)=0$ for $k=1,2,\ldots,n$ \cite[2.24.2.2]{PBM3}:
\begin{equation}\label{eq:GntildeG0}
\widetilde{G}_n(t)=\frac{1}{(n-1)!}\int_t^1G_0(x)(x-t)^{n-1}dx.
\end{equation}

When $n>0$, Definition~1 is motivated by the following argument. Replace $\varphi(t)$ in (\ref{eq:Gphiintegral}) by its Taylor expansion at $t=0$:
$$
\varphi(t)=\sum_{k=0}^{n-1}\frac{\varphi^{(k)}(0)}{k!}t^k+\varphi_n(t),
$$
where $\varphi_n(t)$ is the Taylor remainder. Then assume $\Re(\a),\Re(\psi_p)>0$ and use (\ref{eq:GMellin}) to obtain the right hand side of (\ref{eq:actionG0}), but with the second term replaced by
$$
\frac{\Gamma(\b)}{\Gamma(\a)}\int_0^1G_0(t)\varphi_n(t)dt.
$$
Integrating by parts $n$ times and using $\varphi_n^{(n)}(t)=\varphi^{(n)}(t)$, $\varphi^{(k)}_{n}(0)=\widetilde{G}_{k+1}(1)=0$ for $k=0,1,\ldots,n-1$ and Properties~5 and 6, we obtain (\ref{eq:actionG0}).  Alternatively, use the integral form of the Taylor remainder $\varphi_n(t)$ and exchange the order of integrals.  Therefore, (\ref{eq:Gphiintegral}) equals the right hand side of (\ref{eq:actionG0}) when $\Re(\a)>0$ and $\Re(\psi_p)>0$. Moreover, the right hand side of (\ref{eq:actionG0}) is an analytic function of the parameters $\a$ and meromorphic function of the parameters $\b$ with simple poles at $-b_i\in\N_0$, so that the right hand side of (\ref{eq:actionG0}) gives an expression for the analytic continuation of (\ref{eq:Gphiintegral}) in $\a$ to the domain $\Re(\a)>-n$ and the meromorphic continuation in $\b$ to the domain $\Re(\psi_p)>-n$. Hence, the family of distributions $\G_0=\G_0(\a,\b)$ is analytic in the parameters $\a$ and meromorphic in $\b$  with simple poles at $-b_i\in\N_0$ in the above domain.

\textbf{Remark.}  The regularization defined in (\ref{eq:actionG0}) can be easily seen to equal  the Hadamard finite part of the divergent integral (\ref{eq:Gphiintegral}), see
\cite{CostinFriedman,EstradaKanwal} for details.  However, we observe a new phenomenon here. In general, the Hadamard finite part constructed to overcome divergence at zero does not alter the situation at other points, while formula (\ref{eq:actionG0}) regularizes the integral (\ref{eq:Gphiintegral}) \emph{at both points, $0$ and $1$, simultaneously}.

\begin{theorem}\label{th:G0continous}
$\G_0$  is a continuous linear functional on $\CB^\infty[0,1]$ and its definition is independent of $n$.
\end{theorem}
\textbf{Proof.} Linearity is obvious. For continuity, assume that $\varphi_j\to\varphi$ in $\CB^{\infty}[0,1]$ and estimate
\begin{multline*}
|\langle{\G_0,\varphi_j}\rangle-\langle{\G_0,\varphi}\rangle|=|\langle{\G_0,\varphi_j-\varphi}\rangle|\leq\sum\limits_{k=0}^{n-1}\frac{|(\a)_k|}{|(\b)_k|k!}|\varphi^{(k)}_j(0)-\varphi^{(k)}(0)|
\\
+\max\limits_{x\in[0,1]}|\varphi^{(n)}_j(x)-\varphi^{(n)}(x)|\left\vert\frac{\Gamma(\b)}{\Gamma(\a)}\right\vert\int_0^1\vert\widetilde{G}_n(t)\vert{dt}\to0~\text{as}~j\to\infty
\end{multline*}
by definition of convergence in $\CB^{\infty}[0,1]$, and because the last integral in finite by Properties~5 and 6 in the Appendix.
Finally, write $\G_{0,n}$ for the distribution $\G_0$ with $n$ terms in the sum (\ref{eq:actionG0}) and $\G_{0,m}$ for $m\ne{n}$ terms. By definition we must choose $n,m>-\Re(\a)$.  Assume, without loss of generality, $n>m$ and let $\varphi$ be an arbitrary test function.  Integration by parts yields
\begin{multline*}
\langle{\G_{0,n},\varphi}\rangle-\langle{\G_{0,m},\varphi}\rangle=\sum\limits_{k=m}^{n-1}\frac{(\a)_k}{(\b)_kk!}\varphi^{(k)}(0)
+\frac{\Gamma(\b)}{\Gamma(\a)}\int_0^1\widetilde{G}_n(t)\varphi^{(n)}(t)dt
-\frac{\Gamma(\b)}{\Gamma(\a)}\int_0^1\widetilde{G}_m(t)\varphi^{(m)}(t)dt
\\
=\sum\limits_{k=m}^{n-1}\frac{(\a)_k}{(\b)_kk!}\varphi^{(k)}(0)+\frac{\Gamma(\b)}{\Gamma(\a)}\int_0^1\widetilde{G}_n(t)\varphi^{(n)}(t)dt
\\
+\left.\frac{\Gamma(\b)}{\Gamma(\a)}\widetilde{G}_{m+1}(t)\varphi^{(m)}(t)\right|^1_0
-\frac{\Gamma(\b)}{\Gamma(\a)}\int_0^1\widetilde{G}_{m+1}(t)\varphi^{(m+1)}(t)dt
\\
=\sum\limits_{k=m+1}^{n-1}\frac{(\a)_k}{(\b)_kk!}\varphi^{(k)}(0)+\frac{\Gamma(\b)}{\Gamma(\a)}\int_0^1\widetilde{G}_n(t)\varphi^{(n)}(t)dt
-\frac{\Gamma(\b)}{\Gamma(\a)}\int_0^1\widetilde{G}_{m+1}(t)\varphi^{(m+1)}(t)dt,
\end{multline*}
where we have used $\widetilde{G}_{m+1}'(t)=-\widetilde{G}_{m}(t)$, $\widetilde{G}_{m+1}(1)=0$ by (\ref{eq:GntildeG0}) and
$$
\widetilde{G}_{m+1}(0)=\int\limits_{0}^{1}\widetilde{G}_m(t)dt=\frac{\Gamma(\a+m)}{\Gamma(\b+m)m!}~~\text{(by  Property~7).}
$$
Repeating integration by parts $(n-m)$ times yields $\langle{\G_{0,n},\varphi}\rangle-\langle{\G_{0,m},\varphi}\rangle=0$. \hfill$\square$

The action of the distribution $\G_0$ on the Laplace, generalized Stieltjes and cosine kernel expectedly leads to the generalized hypergeometric functions of the Kummer, Gauss and Bessel type, respectively.

\begin{theorem}\label{th:Gaction}
Suppose complex $\a$ and $\b$, $-\b\notin\N_0$,  satisfy $\Re(\a)>-n$, $\Re(\psi_p)>-n$ for some $n\in\N_0$, where $\psi_p$ is defined in \emph{(\ref{eq:a-psi-defined})}.  Then for all $\sigma\in\C$
\begin{align}
\label{eq:G-Stieltjesaction}
&\langle\G_0,(1+zt)^{-\sigma}\rangle={}_{p+1}F_p\left.\left(\begin{matrix}\sigma,\a\\\b\end{matrix}\right\vert-z\right),
\\
\label{eq:G-Laplaceaction}
&\langle\G_0,e^{-zt}\rangle={}_{p}F_p\left.\left(\begin{matrix}\a\\\b\end{matrix}\right\vert-z\right),
\\
\label{eq:G-cosineaction}
&\langle\G_0,\cos(2\sqrt{zt})\rangle={}_{p-1}F_p\left.\left(\begin{matrix}\a_{[p]}\\\b\end{matrix}\right\vert -z\right),
\end{align}
where in the last formula it is assumed that $a_p=1/2$.  Formulas \emph{(\ref{eq:G-Laplaceaction})}, \emph{(\ref{eq:G-cosineaction})} are valid for all complex $z$, while \emph{(\ref{eq:G-Stieltjesaction})} is true for $z\in\C\!\setminus\!(-\infty,-1]$.
\end{theorem}
\textbf{Proof.} As we explained below Definition~1, $\langle\G_0,\varphi(t)\rangle$ is a representation of the analytic continuation in parameters $\a$ and $\b$.  Further, formulas (\ref{eq:G-Stieltjesaction})-(\ref{eq:G-cosineaction}) are true for $\Re(\a),\Re(\psi_p)>0$ as they reduce to (\ref{eq:Frepr}),(\ref{eq:qFqLaplace}) and (\ref{eq:q-1Fqcosine}), respectively, while the the right hand sides are analytic in $\a$ and $\b$ save the poles. This proves (\ref{eq:G-Stieltjesaction})-(\ref{eq:G-cosineaction}).  Alternatively, an application of $\G_0$ to $\varphi_z(t)=(1+zt)^{-\sigma}$ is immediately seen to lead to the right hand side of (\ref{eq:G-Stieltjes}), while taking $\varphi_z(t)=\exp(-zt)$ yields the right hand side of (\ref{eq:G-Laplace}).  For the Bessel type function we will use that
$$
\cos\left(2\sqrt{zt}\right)={}_0F_1(-;1/2;-zt),~~\text{so that}~\frac{\partial^n}{\partial{t}^n}\cos\left(2\sqrt{zt}\right)=\frac{(-z)^n}{(1/2)_n}{}_0F_1(-;n+1/2;-zt).
$$
Setting $a_p=1/2$, we obtain:
\begin{multline*}
\langle\G_0,\cos(2\sqrt{zt})\rangle
\\
=\sum_{k=0}^{n-1}\frac{(\a)_k(-z)^k}{(\b)_k(1/2)_kk!}+\frac{\Gamma(\b)(-z)^n}{\Gamma(\a)(1/2)_n}
\int_0^1\!\!G^{p+1,0}_{p+1,p+1}\left(t\left|\begin{array}{l}\!\!\b-1+n,n\!\!\\\!\!\a-1+n,0\!\!\end{array}\right.\right){}_0F_1\!\left.\!\left(\begin{matrix}\\n+1/2\end{matrix}\right\vert -zt\!\right)\!dt
\\
=\sum_{k=0}^{n-1}\frac{(\a)_k(-z)^k}{(\b)_k(1/2)_kk!}
+\frac{\Gamma(\b)(-z)^n}{\Gamma(\a)(1/2)_n}\sum\limits_{j=0}^{\infty}\frac{\Gamma(\a+n+j)\Gamma(j+1)(-z)^j}{\Gamma(\b+n+j)\Gamma(n+j+1)j!(n+1/2)_j}
\\
=\sum_{k=0}^{n-1}\frac{(\a)_k(-z)^k}{(\b)_k(1/2)_kk!}
+\sum\limits_{j=0}^{\infty}\frac{(\a)_{n+j}(-z)^{n+j}}{(\b)_{n+j}(1/2)_{n+j}(n+j)!}={}_{p-1}F_p\left.\left(\begin{matrix}\a_{[p]}\\\b\end{matrix}\right\vert -z\right).~~~~~~~~~~~~~~\square
\end{multline*}

\textbf{Remark.} Representations (\ref{eq:G-Stieltjesaction}) and (\ref{eq:G-Laplaceaction}) are, of course, just different ways of writing (\ref{eq:G-Stieltjes}) and (\ref{eq:G-Laplace}), respectively.   Nevertheless, representation (\ref{eq:G-cosine}) is essentially different from (\ref{eq:G-cosineaction}), as seen from the proof.

\subsection{An application: extended Luke's inequalities}

In \cite[Theorem~16]{Luke} Luke found two-sided bounds for the functions ${_pF_p}(\a;\b;x)$ and ${_{p+1}F_p}(\sigma,\a;\b;x)$ under the restrictions $b_i\ge{a_i}>0$, $i=1,2,\ldots,p$. The bounds were presented without proofs, but mentioning that they ''can be easily proved''.  In \cite{KarpJMS} the first author gave two different proofs of Luke's inequalities valid for different sign of $x$ and relaxed the conditions on parameters.  Using the decompositions (\ref{eq:G-Stieltjes}) and (\ref{eq:G-Laplace}) and mimicking the proof from \cite{KarpJMS}, we can extend Luke's inequalities to arbitrary real parameter values. Before formulating the result, let us remind the reader that inequalities like $\a>-3$ and sums like $\a+1$ are understood element-wise and $(\a+n)$ denotes the product $\prod_{j=1}^{p}(a_j+n)$.

\begin{theorem}\label{th:LukeExtended1}
Suppose that  $\a,\b\in\R^p$ are such that $\b$ does not contain non-positive integers. Choose $n\in\N_0$ satisfying $\a,\psi_p\ge-n$ such that $v_{\a',\b'}(t)\ge0$ for $t\in[0,1]$, where  $\a'=(\a+n,1)$, $\b'=(\b+n,n+1)$  and $v_{\a',\b'}$ is defined in \emph{(\ref{eq:vge0})} \emph{(}in particular, it is sufficient that $\b'\!\prec^W\!\a'$\emph{)}.
Then for all real $x$,
\begin{multline}\label{eq:LukeExtended1}
\sum\nolimits_{j=0}^{n-1}\frac{(\a)_jx^{j}}{(\b)_jj!}+(-1)^{\alpha}\frac{(\a)_nx^n}{(\b)_nn!}\exp\!\left(\frac{x(\a+n)}{(n+1)(\b+n)}\right)\leq
{}_{p}F_p\left.\!\!\left(\begin{matrix}\a\\\b\end{matrix}\right\vert x\!\right)
\\[5pt]
\leq\sum\nolimits_{j=0}^{n-1}\frac{(\a)_jx^{j}}{(\b)_jj!}+(-1)^{\alpha}\frac{(\a)_nx^n}{(\b)_nn!}\left(\frac{(e^{x}-1)(\a+n)}{(n+1)(\b+n)}+1\right),
\end{multline}
where
$$
\alpha:=\left\{\!\!\!\begin{array}{l}0,~~\emph{if}~~~\Gamma(\a)\Gamma(\b)x^n\ge0,
\\[3pt]
1,~~\emph{if}~~~\Gamma(\a)\Gamma(\b)x^n<0.
\end{array}\right.
$$
If some elements of $\a$ are non-positive integers then equality holds on both sides of \emph{(\ref{eq:LukeExtended1})} regardless of the value of $\alpha$.
\end{theorem}
\textbf{Proof.} Denote
$$
f_n(x):=\sum\limits_{k=n}^{\infty}\frac{(\a)_kx^k}{(\b)_kk!}={}_{p}F_p\left.\!\!\left(\begin{matrix}\a\\\b\end{matrix}\right\vert x\!\right)-\sum\limits_{k=0}^{n-1}\frac{(\a)_kx^k}{(\b)_kk!}
$$
Assume first that $\psi_p>-n$ and $\a,\b>-n$. Then (\ref{eq:G-Laplace}) can be rewritten as follows:
$$
\frac{\Gamma(\a)}{\Gamma(\b)}\frac{f_n(x)}{x^n}=\int_0^1\phi_x(f(t))\mu(dt),~~\text{where}~~
\mu(dt)=G^{p+1,0}_{p+1,p+1}\left(t\left|\begin{array}{l}\!\!\b+n,n+1\!\!\\\!\!\a+n,1\!\!\end{array}\right.\right)\!\frac{dt}{t},
$$
$\phi_x(t)=e^{xt}$ and $f(t)=t$.  According to Property~9, the condition $v_{\a',\b'}(t)\ge0$  is sufficient for the measure $\mu(dt)$ to be nonnegative. Now we can apply  the integral form of  Jensen's inequality  \cite[Chapter~I, formula (7.15)]{MPF},
\begin{equation}\label{eq:Jensen}
\phi_x\left(\int\nolimits_0^{1}f(t)\mu(dt)\!\Bigg/\!\!\int\nolimits_0^{1}\mu(dt)\right)
\leq\int\nolimits_0^{1}\phi_x(f(t))\mu(dt)\!\Bigg/\!\!\int\nolimits_0^{1}\mu(dt),
\end{equation}
valid for convex $\phi_x$ and  $f$  integrable with respect to a nonnegative measure $\mu$.
Computing
$$
\int\nolimits_0^{1}\mu(dt)=\frac{\Gamma(\a+n)}{\Gamma(\b+n)\Gamma(n+1)},~~~\int\nolimits_0^{1}f(t)\mu(dt)=\frac{\Gamma(\a+n+1)}{\Gamma(\b+n+1)\Gamma(n+2)},
$$
we arrive at
$$
\exp\!\left(\frac{x(\a+n)}{(n+1)(\b+n)}\right)\leq\frac{f_n(x)}{x^n}\frac{(\b)_nn!}{(\a)_n}.
$$
Multiplying this formula by the nonnegative number $(-1)^{\alpha}(\a)_nx^n/[(\b)_nn!]$ and recalling the definition of $f_n(x)$ we obtain the lower bound of (\ref{eq:LukeExtended1}).  Further we apply the converse Jensen's inequality in the form \cite[Theorem~3.37]{PPT}
$$
\int\nolimits_0^{1}\phi_x(f(t))d\mu(t)\!\Bigg/\!\!\int\nolimits_0^{1}\mu(dt)
\leq
(\phi_x(1)-\phi_x(0))\int\nolimits_0^{1}f(t)d\mu(t)\!\Bigg/\!\!\int\nolimits_0^{1}\mu(dt)+1\!\cdot\!\phi_x(0)-0\!\cdot\!\phi_x(1).
$$
Substituting we get
$$
\frac{(\b)_nn!}{(\a)_n}\frac{f_n(x)}{x^n}
\leq\frac{(e^{x}-1)(\a+n)}{(n+1)(\b+n)}+1.
$$
Again, multiplying this formula by the nonnegative number $(-1)^{\alpha}(\a)_nx^n/[(\b)_nn!]$ and recalling the definition of $f_n(x)$ we obtain the uper bound of (\ref{eq:LukeExtended1}).  Finally, if $\psi_p=0$ and/or some of the components of $\a$ are equal to $-n$, the inequality is still true by continuity.  The components of $\b$ cannot be equal to non-positive integers by the hypothesis of the theorem. $\hfill\square$

\smallskip

\begin{theorem}\label{th:LukeExtended2}
Suppose that $\sigma\in\R$, $\a,\b\in\R^p$ are such that $\b$ does not contain non-positive integers. Choose $n\in\N_0$ satisfying $\a,\psi_p\ge-n$ such that $v_{\a',\b'}(t)\ge0$ for $t\in[0,1]$, where  $\a'=(\a+n,1)$, $\b'=(\b+n,n+1)$  and $v_{\a',\b'}$ is defined in \emph{(\ref{eq:vge0})} \emph{(}in particular, it is sufficient that $\b'\!\prec^W\!\a'$\emph{)}.
Then for $x<1$
\begin{multline}\label{eq:LukeExtended2}
\sum\nolimits_{j=0}^{n-1}\frac{(\sigma)_j(\a)_jx^{j}}{(\b)_jj!}+(-1)^{\alpha}\frac{(\sigma)_n(\a)_nx^n}{(\b)_nn!}
\left(1-\frac{x(\a+n)}{(n+1)(\b+n)}\right)^{-\sigma-n}\!\leq\!
{}_{p+1}F_p\left.\!\!\left(\begin{matrix}\sigma,\a\\\b\end{matrix}\right\vert x\!\right)
\\[5pt]
\leq\sum\nolimits_{j=0}^{n-1}\frac{(\sigma)_j(\a)_jx^{j}}{(\b)_jj!}+(-1)^{\alpha}\frac{(\sigma)_n(\a)_nx^n}{(\b)_nn!}
\left(\frac{(1-(1-x)^{\sigma+n})(\a+n)}{(1-x)^{\sigma+n}(n+1)(\b+n)}+1\right),
\end{multline}
where
$$
\alpha:=\left\{\!\!\!\begin{array}{l}0,~~\emph{if}~~~\Gamma(\sigma)\Gamma(\a)\Gamma(\b)x^n\ge0,
\\[3pt]
1,~~\emph{if}~~~\Gamma(\sigma)\Gamma(\a)\Gamma(\b)x^n<0.
\end{array}\right.
$$
If some elements of $\a$ are non-positive integers then equality holds on both sides of \emph{(\ref{eq:LukeExtended2})} regardless of the value of $\alpha$.
\end{theorem}
\textbf{Proof.} Start with representation (\ref{eq:G-Stieltjes}) rewritten as follows:
$$
\frac{\Gamma(\a)f_n(x)}{\Gamma(\b)(\sigma)_nx^n}
=\int_0^1G^{p+1,0}_{p+1,p+1}\left(t\left|\begin{array}{l}\!\!\b+n,n+1\!\!\\\!\!\a+n,1\!\!\end{array}\right.\right)\!\frac{dt}{(1-xt)^{\sigma+n}t},
$$
where
$$
f_n(x)=\sum\limits_{k=n}^{\infty}\frac{(\sigma)_k(\a)_kx^k}{(\b)_kk!}
={}_{p+1}F_p\left.\!\!\left(\begin{matrix}\sigma,\a\\\b\end{matrix}\right\vert x\!\right)-\sum\limits_{k=0}^{n-1}\frac{(\sigma)_k(\a)_kx^k}{(\b)_kk!},
$$
and repeat the steps of the proof of the previous theorem with $\phi_x(t)=(1-x)^{-\sigma-n}$.$\hfill\square$

\paragraph{Acknowledgements.}
We thank Leonid Golinskii for sharing his insights regarding radial positive definite functions. Research of the first author has been supported by the Russian Science Foundation under project 14-11-0002.  The research of the second author has been supported by the Spanish \emph{Ministry of "Econom\'{\i}a y Competitividad"} under  project MTM2014-53178.

\bigskip
\bigskip

\section{Appendix. Definition and properties of Meijer's $G$-function}

Suppose that $0\leq{m}\leq{q}$, $0\leq{n}\leq{p}$ are integers and $\a$, $\b$ are arbitrary
complex vectors, such that $a_i-b_j\notin\N$ for all $i=1,\ldots,n$ and $j=1,\ldots,m$. Meijer's $G$-function
is defined by the Mellin-Barnes integral  of the form
(see \cite[section~12.3]{BealsWong}, \cite[section~5.3]{HTF1}, \cite[chapter~1]{KilSaig}, \cite[section~8.2]{PBM3} or \cite[section~16.17]{NIST})
\begin{equation}\label{eq:G-defined}
G^{m,n}_{p,q}\!\left(\!z~\vline\begin{array}{l}\a\\\b\end{array}\!\!\right)\!\!:=
\\
\frac{1}{2\pi{i}}
\int\limits_{\mathcal{L}}\!\!\frac{\Gamma(b_1\!+\!s)\cdots\Gamma(b_m\!+\!s)\Gamma(1-a_1\!-\!s)\cdots\Gamma(1-a_n\!-\!s)z^{-s}}
{\Gamma(a_{n+1}\!+\!s)\cdots\Gamma(a_p\!+\!s)\Gamma(1-b_{m+1}\!-\!s)\cdots\Gamma(1-b_{q}\!-\!s)}ds,
\end{equation}
where the contour $\L$ is a simple loop that separates the poles of the integrand of the form $b_{jl}=-b_j-l$, $l\in\N_0$ leaving them on the left from the poles of the form $a_{ik}=1-a_i+k$, $k\in\N_0$, leaving them on the right \cite[section~1.1]{KilSaig}. The condition $a_i-b_j\notin\N$  guarantees that the poles are separable.
The contour may have one of the three forms $\L_{-}$, $\L_{+}$ or $\L_{i\gamma}$ described below. Choose any
$$
\varphi_1<\min\{-\Im{b_1},\ldots,-\Im{b_m},\Im(1-a_1),\ldots,\Im(1-a_n)\},
$$
$$
\varphi_2>\max\{-\Im{b_1},\ldots,-\Im{b_m},\Im(1-a_1),\ldots,\Im(1-a_n)\}
$$
and  arbitrary real $\gamma$.    The contour $\L_{-}$ is a left loop lying in the horizontal strip $\varphi_1\leq\Im{s}\leq\varphi_2$. It starts at the point $-\infty+i\varphi_1$, terminates at the point $-\infty+i\varphi_2$ and coincides with the sides of the strip for sufficiently large $|s|$.  Similarly, the contour $\L_{+}$ is a right loop lying in the same strip, starting at the point $+\infty+i\varphi_1$ and terminating at the point $+\infty+i\varphi_2$.  It coincides with the sides of the strip for sufficiently large $|s|$.  Finally, the contour  $\L_{i\gamma}$ starts at $\gamma-i\infty$,  terminates at $\gamma+i\infty$ and coincides with the line $\Re{s}=\gamma$ for all sufficiently large $|s|$.  The power function $z^{-s}$ is defined on the Riemann surface of the logarithm, so that
$$
z^{-s}=\exp(-s\{\log|z|+i\arg(z)\})
$$
and $\arg(z)$ is allowed to take any real value. Hence, $G^{m,n}_{p,q}(z)$ is also defined on the Riemann surface of the logarithm.
Set \cite[(1.1.10)]{KilSaig}
$$
\mu:=\sum\nolimits_{j=1}^{q}b_j-\sum\nolimits_{i=1}^{p}a_i-\frac{p-q}{2}, ~~~~a^*:=2(m+n)-(p+q).
$$
Specialization of \cite[Theorem~1.1]{KilSaig} (which deals with a more general Fox's $H$-functions) to our situation leads to the following conditions for convergence of the integral in (\ref{eq:G-defined}):

(a) if $\L=\L_{-}$ the integral in (\ref{eq:G-defined}) converges for $0<|z|<1$ and arbitrary $\a$, $\b$ and also for $|z|=1$ if $\Re(\mu)<-1$;

(b) if $\L=\L_{+}$ the integral in (\ref{eq:G-defined}) converges  for $|z|>1$ and arbitrary $\a$, $\b$ and also for $|z|=1$ if $\Re(\mu)<-1$;

(c) if $\L=\L_{i\gamma}$ the integral in (\ref{eq:G-defined}) converges for $|\arg(z)|<a^*\pi/2$, $z\ne0$ if $a^*>0$ and $\arg(z)=0$, $z\ne1$, $\Re(\mu)<0$ if $a^*=0$.

\noindent The last condition has been proved in \cite[Theorem~3.3]{KilSaig} in a more general case and earlier in the first author's thesis for $G^{p,0}_{p,p}$ (see also \cite[Lemma~1]{KPJMAA}). Note that \cite[Theorem~1.1]{KilSaig} requires a stronger restriction $\Re(\mu)<-1$.  If the integral in (\ref{eq:G-defined}) exists for several contours, the resulting functions coincide in all known cases. A more detailed discussion of this issue can be found in our recent paper \cite{KPSIGMA}.  A comprehensive overview of the properties of $G$-function is contained in \cite[Section~8.2]{PBM3} and in \cite{Wolfram}.  In this paper we mostly need the properties of $G^{p,0}_{p,p}$ found in the above references as well as some of its new or less obvious properties for which we will supply detailed references or explanations. In what follows we will write $\a_{[k_1,k_2,\ldots,k_r]}$ for the vector $\a$ with the elements $a_{k_1},a_{k_2},\ldots,a_{k_r}$ removed. In particular,  $\a_{[k]}=(a_1,\ldots,a_{k-1},a_{k+1},\ldots,a_p)$.
As before, $\Gamma(\a)$ stands for $\prod_{i=1}^{p}\Gamma(a_i)$ and $\a+\alpha$ with scalar $\alpha$ is an abbreviation for the vector $(a_1+\alpha,\ldots,a_p+\alpha)$.

\textbf{Property~1.} The function $G^{m,n}_{p,q}$ is real if the vectors $\a$, $\b$ and the argument $z$ are real.  This follows from the residue expansion, see \cite[8.2.2.3-4]{PBM3} or \cite[16.17.2]{NIST}.

\textbf{Property~2.} For any real $\alpha$,
$$
z^{\alpha}G^{m,n}_{p,q}\!\left(\!z~\vline\begin{array}{l}\b\\\a\end{array}\!\!\right)=G^{m,n}_{p,q}\!\left(\!z~\vline\begin{array}{l}\b+\alpha\\\a+\alpha\end{array}\!\!\right).
$$
See \cite[8.2.2.15]{PBM3} or \cite[16.19.2]{NIST}. The property also holds for complex $\alpha$, but care must be taken in choosing the correct branches.

\textbf{Property~3.}  According to \cite[Lemma~1]{KPJMAA} and \cite[Theorem~6]{KPCMFT},
$$
G^{p,0}_{p,p}\!\left(\!z~\vline\begin{array}{l}\b\\\a\end{array}\!\!\right)=0~~\text{for}~~|z|>1.
$$
The above $G$-function is well-defined for arbitrary values of $\a$ and $\b$ if the contour $\L$ is chosen to be $\L_{+\infty}$. Under the restriction $\Re(\psi_p)>0$, it can also be deformed into $\L_{i\gamma}$, where $\psi_p=\sum_{j=1}^{p}(b_j-a_j)$.

\textbf{Property~4.} If none of the vectors $\a_{[k]}-a_k$, $k=1,\ldots,p$, contains integers,  Meijer's $G$ function can
be expanded in terms of generalized hypergeometric functions as follows
\begin{equation}\label{eq:Gp0qp}
G^{p,0}_{p,p}\!\left(\!z~\vline\begin{array}{l}\b
\\\a\end{array}\!\!\right)=
\sum\limits_{k=1}^{p}z^{a_k}\frac{\Gamma(\a_{[k]}-a_k)}{\Gamma(\b-a_k)}
{_pF_{p-1}}\!\left(\begin{array}{l}1-\b+a_k
\\1-\a_{[k]}+a_k\end{array}\vline\: z\right).
\end{equation}
See \cite[(34)]{Marichev}, \cite[8.2.2.3]{PBM3} or \cite[16.17.2]{NIST}.

\textbf{Property~5.}  Note that the poles of the numerator of the integrand $z^{-s}\Gamma(\a+s)/\Gamma(\b+s)$ in the definition of $G^{p,0}_{p,p}(z)$ may cancel out with the poles of the denominator.   Suppose that $b_k=a_i+q$ for some $k=1,\ldots,p$ and $q\in\Z$.  If $q\leq0$, then all the poles of the function $\Gamma(a_i+s)$ cancel out with poles of $\Gamma(b_k+s)$. We will call the indices $i$ and the corresponding components of $\a$ normal if at least one pole of $\Gamma(a_i+s)$ does not cancel (if this pole is single then it is the rightmost pole) . We say that $\a$ is normal if all its components are normal. In general situation we can ''normalize'' $\a$ by deleting the exceptional (not normal) components.

Suppose $\a=(a_1,a_2,\ldots,a_{p'})$ is normal or normalized. In general, it  may contain some groups of equal elements (on the extreme all elements are allowed to be equal as well). Write $r$ for the cardinality of the largest group of equal elements for which $\min(\Re(a_1),\ldots,\Re(a_{p'}))$ is attained. Assume for a moment that there is only one such group and suppose, without loss of generality, that this group is $a_1=a_2=\cdots=a_r=a$.   Then
\begin{equation}\label{eq:G-asymp-zero}
G^{p,0}_{p,p}\!\left(z~\vline\begin{array}{l}\!\b\!\\\!\a\!\end{array}\right)={\alpha}z^{a}\log^{r-1}(z)(1+O(\log^{-1}(z)))~~\text{as}~~z\to{0},
\end{equation}
where
$$
\alpha=\frac{(-1)^{r-1}\prod\nolimits_{i=r+1}^{p}\Gamma(a_i-a)}{(r-1)!\prod\nolimits_{i=1}^{p}\Gamma(b_i-a)}.
$$
For $r=1$ the term $O(\log^{-1}(z))$ must be substituted with $O(z^{\delta}\log^k(z))$, where $\delta=\Re(\tilde{a}-a)$ and $\tilde{a}$ is the element with the second smallest real part while $k$ stands for its multiplicity.   If there are several groups of equal elements of the same cardinality $r$ for which $\min(\Re(a_1),\ldots,\Re(a_p'))$ is attained, then formula (\ref{eq:G-asymp-zero}) remains valid with the constant $\alpha$ equal to the sum of the corresponding constants for each group  (computed as above).  Note that $\alpha\ne0$ by normality of $\a$. The asymptotic approximation as $z\to0$ for a more general Fox's $H$ function is given in \cite[Theorem~1.5]{KilSaig}. However, the computation of the constant in \cite[formula~(1.4.6)]{KilSaig} seems to contain an error, corrected in (\ref{eq:G-asymp-zero}) using residue expansion \cite[formula~(1.2.22)]{KilSaig}.

\textbf{Property~6.} An important property used in this paper is the following representation:
\begin{equation}\label{eq:Norlund}
G^{p,0}_{p,p}\!\left(\!z~\vline\begin{array}{l}\b\\\a\end{array}\!\!\right)=\frac{z^{a_k}(1-z)^{\psi_p-1}}{\Gamma(\psi_p)}
\sum\limits_{n=0}^{\infty}\frac{g_n(\a_{[k]};\b)}{(\psi_p)_n}(1-z)^n,~~~k=1,2,\ldots,p,
\end{equation}
which holds in the disk $|1-z|<1$ for all $-\psi_p=-\sum_{i=1}^{p}(b_i-a_i)\notin\N_0$ and each $k=1,2,\ldots,p$. Several ways are known to compute the coefficients $g_n(\a_{[k]};\b)$. They satisfy two different recurrence relations (in $p$ and $n$). The simplest of them reads
\begin{equation}\label{eq:Norlundcoeff}
g_n(\a_{[p+1]};\b)=\sum\limits_{s=0}^{n}\frac{(b_{p+1}-a_{p})_{n-s}}{(n-s)!}(\psi_{p}+s)g_s(\a_{[p,p+1]};\b_{[p+1]}),~~~p=1,2,\ldots,
\end{equation}
with initial values $g_0(-;b_1)=1$, $g_n(-;b_1)=0$, $n\ge1$.
The coefficient $g_n(\a_{[k]};\b)$ is obtained from $g_n(\a_{[p]};\b)$ by
exchanging the roles of $a_p$ and $a_k$, or by using the connection formula
\begin{equation}\label{eq:Norlundconnection}
g_n(\a_{[k]};\b)=\sum\limits_{s=0}^{n}\frac{(a_k-a_p)_{n-s}}{(n-s)!}(\psi_{p}+s)g_s(\a_{[p]};\b),~~~k=1,2,\ldots,p.
\end{equation}
The following explicit representation was derived in \cite[(1.28), (2.7), (2.11)]{Norlund}:
\begin{equation}\label{eq:Norlund-explicit}
g_n(\a_{[p]};\b)=\sum\limits_{0\leq{j_{1}}\leq{j_{2}}\leq\cdots\leq{j_{p-2}}\leq{n}}
\prod\limits_{m=1}^{p-1}\frac{(\psi_m+j_{m-1})_{j_{m}-j_{m-1}}}{(j_{m}-j_{m-1})!}(b_{m+1}-a_{m})_{j_{m}-j_{m-1}},
\end{equation}
where $\psi_m=\sum_{i=1}^{m}(b_i-a_i)$, $j_0=0$, $j_{p-1}=n$.
Expansion (\ref{eq:Norlund}) using different notation and without mentioning $G$-function was derived by N{\o}rlund in \cite[formulas~(1.33), (1.35), (2.7)]{Norlund}.
The history and many further details regarding N{\o}rlund's results and methods to compute $g_n(\a_{[k]};\b)$ can be found in our recent paper \cite{KPSIGMA}.

Taking limit $\psi_p\to-l$, $l\in\N_0$ in (\ref{eq:Norlund}) we obtain
\begin{equation}\label{eq:Norlund1}
G^{p,0}_{p,p}\!\left(\!z~\vline\begin{array}{l}\b\\\a\end{array}\!\!\right)=z^{a_k}
\sum\limits_{n=0}^{\infty}\frac{g_{n+l+1}(\a_{[k]};\b)}{n!}(1-z)^n,~~~k=1,2,\ldots,p,
\end{equation}
where $\psi_p=-l$, $l\in\N_0$ (see \cite[formula~(1.34)]{Norlund}). Hence, $G^{p,0}_{p,p}$ is analytic in the neighborhood of $z=1$ for non-positive integer values of $\psi_p$.

\textbf{Property~7.}  The Mellin transform of $G^{p,0}_{p,p}$ exists if either  $\Re(\psi_p)>0$ or
$\psi_p=-m$, $m\in\N_0$. In the former case
\begin{equation}\label{eq:GMellin}
\int\limits_{0}^{\infty}x^{s-1}G^{p,0}_{p,p}\!\left(\!x~\vline\begin{array}{l}\b\\\a\end{array}\!\!\right)dx
=\int\limits_{0}^{1}x^{s-1}G^{p,0}_{p,p}\!\left(\!x~\vline\begin{array}{l}\b\\\a\end{array}\!\!\right)dx
=\frac{\Gamma(\a+s)}{\Gamma(\b+s)}
\end{equation}
is valid in the intersection of the half-planes $\Re(s+a_i)>0$ for $i=1,\ldots,p$.  If $\psi_p=-m$, $m\in\N_0$ then
\begin{equation}\label{eq:GMellinNorlund}
\int\limits_{0}^{\infty}x^{s-1}G^{p,0}_{p,p}\!\left(\!x~\vline\begin{array}{l}\b\\\a\end{array}\!\!\right)dx
=\int\limits_{0}^{1}x^{s-1}G^{p,0}_{p,p}\!\left(\!x~\vline\begin{array}{l}\b\\\a\end{array}\!\!\right)dx
=\frac{\Gamma(\a+s)}{\Gamma(\b+s)}-q(s)
\end{equation}
in the same half-plane.  Here $q(s)$ is a polynomial of degree $m$ given by
\begin{equation}\label{eq:q-polynomial}
q(s)=\sum\limits_{j=0}^{m}g_{m-j}(\a_{[k]};\b)(s+a_k-j)_j,~~~k=1,2,\ldots,p.
\end{equation}
The coefficients $g_{i}(\a_{[k]};\b)$ depend on $k$.  The resulting polynomial $q(s)$, however, is the same for each $k$.  See \cite[(2.18), (2.29)]{Norlund} or \cite[(4)]{KPFOX}.

\textbf{Property~8.} Given a nonnegative integer $k$ suppose that $\Re(\psi_p)>-k$ and $\Re(a_i)>0$ for $i=1,\ldots,p$. Then we have
\begin{equation}\label{eq:intG1-x}
\int\limits_{0}^{1}G^{p,0}_{p,p}\!\left(\!x~\vline\begin{array}{l}\b-1\\\a-1\end{array}\!\!\right)(1-x)^kdx
=\frac{\Gamma(\a)}{\Gamma(\b)}{}_{p+1}F_p\left.\left(\begin{matrix}-k, \a\\ \b\end{matrix}\right\vert 1\right).
\end{equation}
Formulas (\ref{eq:G-asymp-zero}) and (\ref{eq:Norlund}) confirm that the integral converges for the specified range of parameters. To demonstrate the validity of (\ref{eq:intG1-x}) we assume first that  $\Re(\psi_p)>0$.  Then the binomial expansion of $(1-x)^k$ and an application of (\ref{eq:GMellin}) yields (\ref{eq:intG1-x}).  Analytic continuation in $\psi_p$ extends the formula to $\Re(\psi_p)>-k$.  We emphasize  that unlike Property~7 formula (\ref{eq:intG1-x}) remains true for non-positive integer $\psi_p>-k$. This can also be seen directly from (\ref{eq:GMellinNorlund}). Indeed, combination of this formula with  the binomial expansion of $(1-x)^k$ gives:
\begin{equation}\label{eq:G(1-x)^k}
\int\limits_{0}^{1}G^{p,0}_{p,p}\!\left(\!x~\vline\begin{array}{l}\b-1\\\a-1\end{array}\!\!\right)(1-x)^kdx
=\frac{\Gamma(\a)}{\Gamma(\b)}{}_{p+1}F_p\left.\left(\begin{matrix}-k, \a\\ \b\end{matrix}\right\vert 1\right)
-\sum\limits_{j=0}^{k}(-1)^j\binom{k}{j}q(j).
\end{equation}
Here
$$
\sum\limits_{j=0}^{k}(-1)^j\binom{k}{j}q(j)=\Delta^kq(0),
$$
where $\Delta{q(s)}:=q(s+1)-q(s)$, $\Delta^kq(s):=\Delta(\Delta^{k-1}q(s))$. But since $q(s)$ has degree $-\psi_p<k$, then $\Delta^kq(0)=0$ confirming (\ref{eq:intG1-x}).  Note that an analogous formula holds for non-integer $\lambda>-\Re{\psi_p}$:
$$
\int\limits_{0}^{1}G^{p,0}_{p,p}\!\left(\!x~\vline\begin{array}{l}\b-1\\\a-1\end{array}\!\!\right)(1-x)^{\lambda}dx
=\frac{\Gamma(\a)}{\Gamma(\b)}{}_{p+1}F_p\left.\left(\begin{matrix}-\lambda, \a\\ \b\end{matrix}\right\vert 1\right).
$$
Indeed, condition $\lambda>-\Re{\psi_p}$ guarantees the convergence of the series on the right \cite[Section 16.2(iii)]{NIST}, so that the formula holds for $\Re{\psi_p}>0$ by the binomial theorem and termwise integration and for $\Re{\psi_p}>-\lambda$ by analytic continuation.  Direct verification of the case $\psi_p=-m$ leads to the identity
$$
\int\limits_{0}^{1}G^{p,0}_{p,p}\!\left(\!x~\vline\begin{array}{l}\b-1\\\a-1\end{array}\!\!\right)(1-x)^{\lambda}dx
=\frac{\Gamma(\a)}{\Gamma(\b)}{}_{p+1}F_p\left.\left(\begin{matrix}-\lambda, \a\\ \b\end{matrix}\right\vert 1\right)
-\sum\limits_{j=0}^{\infty}(-1)^j\binom{\lambda}{j}q(j).
$$
But  for each $m<\lambda$,
$$
\sum\limits_{j=0}^{\infty}(-1)^j\binom{\lambda}{j}j^m=0.
$$

\textbf{Property~9.} As we mentioned in the introduction, the inequality
$$
G^{p,0}_{p,p}\!\left(\!x~\vline\begin{array}{l}\b\\\a\end{array}\!\!\right)\ge0~~\text{for}~~0<x<1
$$
holds if $v_{\a,\b}(t)=\sum_{j=1}^{p}(t^{a_j}-t^{b_j})\ge0$ for $t\in[0,1]$. See \cite[Theorem~2]{KarpJMS} for a proof of this fact and \cite[section~2]{KPCMFT} for further details.  Note also that $v_{\a,\b}(t)\ge0$ implies that $\psi_p=\sum_{j=1}^{p}(b_j-a_j)\ge0$.  For given $\a$, $\b$ the inequality $v_{\a,\b}(t)\ge0$ is not easy to verify other than numerically. However, several sufficient conditions for $v_{\a,\b}(t)\ge0$ expressed directly in terms of $\a$, $\b$ are known.   In particular, according to  \cite[Theorem~10]{Alzer} $v_{\a,\b}(t)\ge0$ on $[0,1]$ if
\begin{equation}\label{eq:amajorb}
\begin{split}
& 0<a_1\leq{a_2}\leq\cdots\leq{a_p},~~
0<b_1\leq{b_2}\leq\cdots\leq{b_p},
\\
&\sum\limits_{i=1}^{k}a_i\leq\sum\limits_{i=1}^{k}b_i~~\text{for}~~k=1,2\ldots,p.
\end{split}
\end{equation}
These inequalities are known as weak supermajorization \cite[Definition~A.2]{MOA} and are abbreviated as
$\b\!\prec^W\!\a$.  Different conditions have been found in \cite[Theorems~1.1,1.2]{GrinIsm}. We slightly generalized the results of  \cite{GrinIsm} and made a survey of other conditions sufficient for $v_{\a,\b}(t)\ge0$ on $[0,1]$ in \cite[section~2]{KPCMFT}. In particular, \cite[Theorem~1]{KPCMFT} asserts that $v_{\a,\b}(t)\ge0$ on $[0,1]$ if $p=2^{n-1}$
and
$$
\a=\biggl\{\sum_{i\in{J}}\alpha_i+\sum_{i\in{I_n\setminus{J}}}\beta_i:~\text{for all}~J\subset I_n=\{1,2,\ldots,n\}~\text{containing even number of terms}\biggr\},
$$
$$
\b=\biggl\{\sum_{i\in{J}}\alpha_i+\sum_{i\in{I_n\setminus{J}}}\beta_i:~\text{for all}~J\subset I_n=\{1,2,\ldots,n\}~\text{containing odd number of terms}\biggr\},
$$
where $\alpha_i\ge\beta_i\ge0$ for $i=1,\ldots,n$.  For example, for $n=3$:
$$
\a=\biggl(\beta_1+\beta_2+\beta_3, \beta_1+\alpha_2+\alpha_3, \alpha_1+\beta_2+\alpha_3, \alpha_1+\alpha_2+\beta_3\biggr),
$$
$$
\b=\biggl(\alpha_1+\beta_2+\beta_3,\beta_1+\alpha_2+\beta_3,\beta_1+\beta_2+\alpha_3,\alpha_1+\alpha_2+\alpha_3\biggr).
$$
Furthermore, these conditions may be combined with (\ref{eq:amajorb}), i.e. if $\a$, $\b$ are given in the above example while $\b_1\prec^W\a_1$ then $v_{\a',\b'}(t)\ge0$ on $[0,1]$ for $\a'=(\a,\a_1)$, $\b'=(\b,\b_1)$.

\end{document}